\documentclass[11pt,reqno,a4paper]{amsart}

\usepackage{amsfonts}
\usepackage{epsfig}
\usepackage{mathrsfs}
\usepackage{amsmath}
\usepackage{amssymb}
\usepackage{color}
\usepackage{mathrsfs}
\usepackage{tikz}
\usepackage{soul}
\usepackage{txfonts}
\usepackage{xparse}
\usepackage[colorlinks=true,linkcolor=blue, urlcolor=red, citecolor=blue]{hyperref}
\usepackage{xcolor}
\usepackage{graphicx}
\usepackage{subcaption}
\usepackage{chngcntr}
\usepackage[T1]{fontenc}
\newtheorem{theorem}{Theorem}[section]
\newtheorem{proposition}[theorem]{Proposition}
\newtheorem{lemma}[theorem]{Lemma}
\newtheorem{corollary}[theorem]{Corollary}
\newtheorem{example}[theorem]{Example}

\newtheorem{remark}[theorem]{Remark}

\newtheorem{definition}[theorem]{Definition}
\newtheorem{maintheorem}{Theorem}

\newcommand{\R}{\mathbb{R}}

\newcommand{\w}{\omega}

\DeclareMathAlphabet{\mathpzc}{OT1}{pzc}{m}{it}

\newcommand{\comment}[1]{}

\usepackage[normalem]{ulem}

\title[ On the stability of $\ddot x(t)+\alpha(t)\dot x(t)+\beta(t) x(t)=0$]
{On the stability of $\ddot x(t)+\alpha(t)\dot x(t)+\beta(t) x(t)=0$}

\author[M. Bessa]{M\'{a}rio Bessa}
\address{Universidade Aberta - Departamento de Ci\^encias e Tecnologia
\\ \newline
Rua do Amial 752, 4200-055 Porto, Portugal}
\email{mario.costa@uab.pt}

\author[H. Vilarinho]{Helder Vilarinho}
\address{Centro de Matem\'atica e Aplica\c{c}\~oes (CMA-UBI), Universidade da Beira
	Interior, Rua Marqu\^es d'\'Avila e Bolama, 6201-001, Covilh\~a, Portugal}
\email{helder@ubi.pt}

\date{\today}

\begin{document}
	
	\begin{abstract}
	Our main goal is to understand the stability of second order linear homogeneous differential equations $\ddot x(t)+\alpha(t)\dot x(t)+\beta(t)x(t)=0$ for $C^0$-generic values of the variable parameters $\alpha(t)$ and $\beta(t)$. For that we embed the problem into the framework of the general theory of continuous-time linear cocycles induced by the random ODE $\ddot x(t)+\alpha(\varphi^t(\omega))\dot x(t)+\beta(\varphi^t(\omega))x(t)=0$, where the coefficients $\alpha$ and $\beta$ evolve along the $\varphi^t$-orbit for $\omega\in M$, and $\varphi^t: M\to M$ is a flow defined on a compact Hausdorff space $M$ preserving a probability measure $\mu$. Considering $y=\dot x$, the above random ODE can be rewritten as $\dot X=A(\varphi^t (\w))X$, with $X=(x,y)^\top$, having a kinetic linear cocycle as fundamental solution.	
	We prove that for a $C^0$-generic choice of parameters $\alpha$ and $\beta$ and for $\mu$-almost all $\omega\in M$ either the Lyapunov exponents of the linear cocycle are equal ($\lambda_1(\omega)=\lambda_2(\omega)$), or else the orbit of $\omega$ displays a dominated splitting. Applying to dissipative systems ($\alpha<0$) we obtain a dichotomy: either $\lambda_1(\omega)=\lambda_2(\omega)<0$, attesting the stability of the solution of the random ODE above, or else the orbit of $\omega$ displays a dominated splitting. Applying to frictionless systems ($\alpha=0$) we obtain a dichotomy: either $\lambda_1(\omega)=\lambda_2(\omega)=0$, attesting the asymptotic neutrality of the solution of the random ODE above, or else the orbit of $\omega$ displays a hyperbolic splitting attesting the \emph{uniform} instability of the solution of the ODE above. This last result implies also an analog result for the 1-d continuous aperiodic Schr\"odinger equation. Furthermore, all results hold for $L^\infty$-generic parameters $\alpha$ and $\beta$.
	\end{abstract}

	\maketitle
	
	\noindent
	
	
	\noindent 
	
	\bigskip
	
	\textbf{\textup{2010} Mathematics Subject Classification:}
	
	Primary: 34D08,  37H15,
	Secondary: 34A30, 37A20.

	\section{Introduction, some definitions and statement of the results}

	\subsection{Introduction} Second-order linear homogeneous differential equations, such as 
	\begin{equation}\label{SOLH}
\ddot x(t)+\alpha(t)\dot x(t)+\beta(t) x(t)=0,
\end{equation}
where $x(t),\alpha(t)$ and $\beta(t)$ are real functions are widely used in physics, engineering, biology, and va\-rious other fields of mathematics. These equations have numerous applications, and some celebrated examples include Hill's, Mathieu's, Meissner's, Lam\'e's and Heun's differential equations and also the time-independent 1-d Schr\"odinger's differential equation. One typical example of their application is in the dynamics of a simple pendulum, which consists of a ball of mass $m$ suspended by a weightless rigid string of length $\ell$ from a fixed support point and subject to gravitational forces. Two noteworthy cases of the simple pendulum include:
\begin{itemize}
\item the \emph{frictionless case} which is described by 
$$\ddot x(t)+\frac{g}{\ell} \sin x(t)=0,$$ where $g$ stands for the gravitational constant, and can be reduced to an equation like \eqref{SOLH} by considering tiny oscillation amplitudes where the simplification $\sin x(t)\sim x(t)$ makes sense, and
\item the \emph{damped case} which is described by 
$$\ddot x(t)+\frac{b}{m} \dot x(t)+\frac{g}{\ell} x(t)=0,$$ where the damping $b$ comes coupled with a first order term which is a typical as e.g. the parachute problem where the damping arises from a first order term related to air resistance. 
\end{itemize}
It is well-known that the general solution of \eqref{SOLH} when $\alpha$ and $\beta$ are constants  depends on the roots of the characteristic corresponding quadratic equation. Despite the problem increases dramatically its difficulty when allowing variation of the parameters $\alpha$ and $\beta$ we still have equations of the form \eqref{SOLH}, e.g. the Cauchy-Euler ones, that can be solved by quadratures.

From the purely existential point of view we have the initial-value problem for \eqref{SOLH} which consists on finding a solution $x(t)$ of the differential equation that also satisfies initial conditions $x(t_0)=x_0$ and $\dot x(t_0)=x_1$. This problem has an affirmative solution if $\alpha(t)$ and $\beta(t)$ are continuous on a certain interval. In between these two approches, say: (a) the exhibition of an analytic solution $x(t)$ and (b) the (unsatisfactory) certainty that the solution exists, relies a fruitful qualitative asymptotic analysis on the behaviour of a solution $x(t)$ although we have no idea of the explicit expression of $x(t)$.

\subsection{Setting up the scenario}

In the present work we intend to decribe with a certain degree of accuracy the asymptotic behavior of $x(t)$, the solution of \eqref{SOLH}, for a generic subset of choices of the parameters $\alpha(t)$ and $\beta(t)$. That is, taking as example the aformentioned simple  pendulum we will be able to describe the limit dynamics of $x(t)$ not for the variable mass $m(t)$, or length $\ell(t)$, or gravity $g(t)$ or even friction $b(t)$ but for respectivelly arbitrarilly close choices $\tilde m(t)$, $\tilde\ell(t)$, $\tilde g(t)$ or $\tilde b(t)$ where \emph{close} means the uniform convergence norm. We follow the steps of the discipline of qualitative theory of differential equations created by Poincar\'e and Lyapunov which pops up as an alternative to the feeble approach of applying analytic methods to integrate most functions confirmed by Liouville's theory. We rewrite \eqref{SOLH} as
	\begin{equation}\label{SOLH2}
		\ddot{x}(t) +\alpha(\varphi^t(\omega))\dot{x}(t)+\beta(\varphi^t(\omega))x(t)=0,
	\end{equation}
	where $\varphi^t\colon M\rightarrow M$ for $t\in\mathbb{R}$ and $\w\in M$ stands for a given flow, say a $\mathbb{R}$-action. This allows, instead of deal with a single equation, to consider infinite equations simultaneously each one for each orbit $\varphi^t(\w)$. The qualitative analysis will be on the Lyapunov exponents of the matricial solution $U(\w,t)$ associated to the linear variational equation 
		\begin{equation}\label{LVE}
		\dot{U}(\w,t)= A(\varphi^{t}(\omega))\cdot U(\w,t),
	\end{equation}
with infinitesimal generator 
\begin{equation}\label{LVE2}
		A(\omega)=\begin{pmatrix}0 & 1\\ -\beta(\omega) & -\alpha(\omega)\end{pmatrix}.
	\end{equation}
Clearly, when $\alpha$ and $\beta$ are periodic coeffficients Floquet theory help us in the analysis and when $\alpha$ and $\beta$ are first integrals, i.e. are constant along the orbits of the flow $\varphi^t$, then \eqref{SOLH2} can be solved by elementary algorithms present in any differential equations book. The interesting case here is when the parameters vary in time along nonperiodic orbits.

		\subsection{Yet another Ma\~n\'e-Bochi dichotomy}
	
	Our objective is to understand the asymptotic behavior of the solutions of \eqref{LVE} when allowing a $C^0$-small perturbation on the parameters of \eqref{SOLH2}. There are several contributions on the literature with respect to this problem \cite{Nov,Mi, F, F2, FJ1, FJ2, FJZ, Be, Be2} sometimes generalizing our case, sometimes considering more restrictive contexts. Two results deserve special attention, as they are closely related to our work. In \cite{Be} the first author inspired in the Ma\~n\'e-Bochi dichotomy (see \cite{M1,M2,B}) proved that $C^0$-generically two-dimensional traceless linear differential systems over a conservative flow have, for almost every point, a dominated splitting or else a trivial Lyapunov spectrum. In \cite{F} Fabbri proved that Schr\"odinger cocycles\footnote{In \eqref{LVE2} take $\alpha=0$ and $\beta=-E+Q(\varphi^t(\w))$ where $E$ is the energy and $Q$ a quasi-periodic potential. See \S\ref{Sc} for more details.} with a quasi-periodic potential and over a certain flow on the torus display a similar dichotomy. These results are very close to our frictionless case setting. However, in \cite{Be} the linear differential system evolve in a broader family, and in \cite{F} the linear differential systems evolve as we saw in a much more rigid family. We intend to pursue a result in the line of Ma\~n\'e-Bochi dichotomy. Clearly, the perturbation framework developed in \cite{Be} can not be used here. Indeed, the fact that the orthogonal group is contained in the special linear group $\text{SL}(2,\mathbb{R})$ was crucial since most perturbations were basically rotating Oseledets directions. Rotating is a much more delicate issue as we will see. Furthermore, the two-dimensional assumption play a crucial role because it seems there is no hope to develop analog results in our setting, like the ones in \cite{BV2,Be2} for example, since apparently there is no chance to mimic a rotational behaviour in higher dimensions. 
	
\subsection{Statement of the main results and some ideas underlying the proofs}

	Notice that infinitesimal generators like $A$ in \eqref{LVE2} gives rise to a particular class of solutions. Clearly, when $\alpha\not=0$ the solutions of \eqref{LVE} evolve on a subclass of the general linear group $\text{GL}(2,\mathbb{R})$ and when $\alpha=0$ the solutions evolve on a subclass of $\text{SL}(2,\mathbb{R})$. Yet both subclasses are not subgroups. Therefore, a particular study must be done taking into consideration that perturbations \textbf{must} belong to our class and not to the broader class of cocycles evolving in $\text{GL}(2,\mathbb{R})$ or even in $\text{SL}(2,\mathbb{R})$. Questions related to this particular class were treated in several works like e.g. \cite{ACE, AW, Be3, FL, Lei,ABV,ABV2}. We will describe the Lyapunov spectrum taking into account the possibility of making a $C^0$-type perturbation on its coefficients. We can resume our initial setup as:
	\begin{itemize}
	\item In the base of the continuous-time cocycle we will consider a $\mathbb{R}$-action $\varphi$ on a compact Hausdorff space $M$ leaving invariant  a Borel regular probability measure $\mu$.
	\item In the fiber of the cocycle we consider the infinitesimal generator of the form \eqref{LVE2} where $\alpha$ and $\beta$ are $C^0$ functions.
	\item We endow $A\colon M\rightarrow \mathpzc{K}$ with the $C^0$ norm, where $\mathpzc{K}$ are the $2\times2$ matrices of the form \eqref{LVE2}, and let $C^0(M,\mathpzc{K})$ be the set of $C^0$ infinitesimal generators of the form \eqref{LVE2}.
		\end{itemize}

	From a kinetic point of view given the \emph{position} and the \emph{momentum} $(x(0),\dot x(0))$ we intend to study the asymptotic behavior when $t\rightarrow\infty$ of the pair $(x(t),\dot x(t))$ namely asymptotic exponential growth rate given by the \emph{Lyapunov exponent}. In the present work we intend to answer the  question of knowing if is it possible to perturb the coefficients $\alpha$ and $\beta$, in the uniform norm $C^0$, in order to obtain equal Lyapunov exponents. We begin be considering the damped case thus encompassing a wide class of linear second order homogeneous differential equations. The definition of uniform hyperbolicity and dominated splitting we refer to the inequalities \eqref{decdom} and \eqref{hypdecdom}. Next result gives a positive answer to open question 1 in \cite{Be4}.

\begin{maintheorem}\label{tps}
		Let $\varphi^t:M\to M$ be a flow preserving the measure $\mu$. There exists a $C^0$ residual set $\mathcal{R}\subset C^0(M,\mathpzc{K})$, such that if $A\in \mathcal{R}$, then for $\mu$-a.e. $\omega$ either 
		\begin{itemize}
		\item there exists a single Lyapunov exponent or else
		\item the splitting along the orbit of $\w$ is dominated.
		\end{itemize}
	\end{maintheorem}
	
	Then we consider the frictionless case answering to open question 2 in \cite{Be4}:
	
		\begin{maintheorem}\label{tps2}
		Let $\varphi^t:M\to M$ be a flow preserving the measure $\mu$. There exists a $C^0$ residual set $\mathcal{R}\subset C^0(M,\mathpzc{K}_{\,\,\Delta=0})$, such that if $A\in \mathcal{R}$, then either 
		\begin{itemize}
		\item  all Lyapunov exponents vanishes or else
		\item the splitting along the orbit of $\w$ is uniformly hyperbolic.
		\end{itemize}
	\end{maintheorem}

	Finally, we consider the dissipative case following a Lyapunov stability approach:
		\begin{maintheorem}\label{tps3}
		Let $\varphi^t:M\to M$ be a flow preserving the measure $\mu$. There exists a $C^0$ residual set $\mathcal{R}\subset C^0(M,\mathpzc{K}_{\,\,\Delta<0})$, such that if $A\in \mathcal{R}$, then either 
		\begin{itemize}
		\item the solution of almost every equation \eqref{SOLH} is stable existing a single negative Lyapunov exponent or
		\item the solution of every associated equation \eqref{SOLH} is stable existing a dominated splitting with two negative Lyapunov exponents or else
		\item for an open and dense set of initial conditions the solution of every associated equation \eqref{SOLH} is (uniformly) unstable.
		\end{itemize}
	\end{maintheorem}

Notice that our results provide a residual subset which, by Proposition~\ref{complete}, is a dense $G_\delta$. 

		The central idea of the proofs of previous theorems is conceptually contained in the classical approach on \cite{M1, M2,B,Be}. More precisely, that in the absence of a dominated splitting we can go perturbing via small rotations in order to be able to send the direction of the top Lyapunov exponent into the direction of the other Lyapunov exponent. Even though these two directions are very different when all these small angle rotations are added up it is possible to send one into the other. This mixture of directions causes an equilibrium of the rates given by the Lyapunov exponents resulting in a trivial spectrum. More precisely and taking as example Theorem~\ref{tps} we first consider that we have a continuity point of a map that sends $A$ into an integral of the top Lyapunov exponent of $A$ in a region without a dominated splitting. Then we show that such point must have trivial Lyapunov exponent because otherwise, and with an arbitrarilly small $C^0$ perturbation of $A$, we would obtain $B$ performing previous mentioned mixing of directions argument. Let's briefly summarize proof's scheme:

\begin{itemize}
 \item By Proposition ~\ref{access} kinetic cocycles are accessible.
\item By\footnote{Similar results for the time-continuous case are available in the literature, see e.g. \cite{Be,Be2}.} \cite[Theorem 5]{BV2} we know that accessible cocycles are points of continuity of $LE$ defined in \eqref{uscf} iff the Oseledets splitting of the cocycle at $x$ is either dominated or trivial at $\mu$-a.e. $x\in M$.
\item  Lemma~\ref{usc} gives that the function defined in \eqref{uscf} is upper semicontinuous.
\item Theorem~\ref{tps} follows from the fact that the continuity points of an upper semicontinuous function is a residual subset.
\item Theorem~\ref{tps2} follows from the fact that hyperbolicity is equivalent to dominated splitting in 2-dim and trivial spectrum is null spectrum for conservative systems.
\item Theorem~\ref{tps3} follows from Theorem~\ref{tps} and a reinterpretation of dominated splitting under dissipative hypothesis.
\end{itemize}
	
		\section{Basic definitions}

	\subsection{Kinetic linear cocycles}

	Let $M$ be a compact Hausdorff space, $\mu$ a probability measure of $M$, $\mathcal{M}$ the Borel $\sigma$-algebra $M$ and let $ \varphi\colon \R \times M \to M$ be a flow, i.e. an $\mathbb{R}$-action, in the sense that it is a measurable map and
	\begin{itemize}
		\item $\varphi^t \colon M \to M$ given by $\varphi^t (\w) = \varphi(t,\w)$ preserves the measure $\mu$ for all $t \in \R$ and all $\w\in M$ and
		
		\item $\varphi^0 (\w)= \w$ and $\varphi^{t+s}=\varphi^t\circ\varphi^s$ for all $t,s \in\R$ and $\w\in M$.
		\end{itemize}

Let $\mathcal{B}(X)$ be the Borel $\sigma$-algebra of a topological space $X$. A continuous-time linear \emph{random dynamical system} on $(\R^2,\mathcal B(\R^2))$, or a continuous-time \emph{linear cocycle}, over $\varphi$ is a $(\mathcal{B}(\R)\otimes \mathcal{M}/\mathcal{B}(\text{GL}(2, \mathbb{R}))$-measurable map 	
	\begin{equation*}
	\Phi:\R\times M \to \text{GL}(2, \mathbb{R})
	\end{equation*}
	such that the mappings $\Phi(t,\w)$ forms a cocycle over $\varphi$, i.e.,
		\begin{enumerate}
			\item $\Phi(0,\w)=\text{Id}$ for all $\w\in M$;
			\item $\Phi(t+s,\w)=\Phi(t,{\varphi^{s}(\w)})\circ\Phi(s,\w)$, for all $s, t\in\R$ and $\w\in M$,
		\end{enumerate}
		and $t\mapsto \Phi(t,\w)$ is continuous for all $\w\in M$. We recall that having $\w\mapsto\Phi(t,\w)$ measurable for each $t\in\R$ and $t\mapsto \Phi(t,\w)$ continuous for all $\w\in M$ implies that $\Phi$ is measurable in the product measure space. This can be stated for more general settings, but this is enough for our purposes. We present below the two most important classes of examples that we will deal with throughout this paper.

Let $A\colon M\rightarrow \mathbb{R}^{2\times 2}$ be a measurable map, where $\mathbb{R}^{2\times 2}$ stands for the set of $2\times 2$ matrices with real entries, and $\Phi_A$ a cocycle satisfying
\begin{equation}\label{eq:LDS0}
		\Phi_{A}(t,\omega)=\text{Id}+\int_{0}^{t}A(\varphi^{s}(\omega))\Phi_{A}(s,\omega)\,ds.
	\end{equation}
The map $\Phi_A(t,\w)$ is called \emph{the Carath\'eodory solution} or \emph{weak solution} of the matricial linear variational equation
	\begin{equation}\label{LVE10}
		\dot{U}(\w,t)= A(\varphi^{t}(\omega))\cdot U(\w,t).
	\end{equation}
 If the solution $\Phi_A(t,\w)$ is differentiable in time (i.e. with respect to $t$) and satisfies for all $t$ and $v\in\mathbb{R}^2$
	\begin{equation}\label{eq:LDS2}
		\frac{d}{dt}\Phi_A(t,\w)v=A(\varphi^t(\w))\Phi_A(t,\w)v\,\,\,\,\,\,\,\,\,\,\text{and}\,\,\,\,\,\,\,\,\,\, \Phi_A(0,\w) v=v,
	\end{equation}
	then it is called a \emph{classical solution}.
	Of course that $t\mapsto \Phi_A(t,\w)v\,$ is continuous for all $\w$ and $v$. Due to \eqref{eq:LDS2} we call $A:M\to \mathbb{R}^{2\times 2}$ an \emph{infinitesimal generator} of $\Phi_A$. Sometimes, due to the relation between $A$ and $\Phi_A$, we refer to both $A$ and $\Phi_A$ as a linear cocyle/RDS/linear differential system.

\begin{example} \textbf{Continuous linear differential systems:} Let $C^0(M,\mathbb{R}^{2\times 2})$ be the set of  $C^0$ maps $A\colon M\rightarrow \mathbb{R}^{2\times 2}$. For a given  $A\in C^0(M,\mathbb{R}^{2\times 2})$ equation~\eqref{LVE10} generates a classical solution $\Phi_A$  which is a $C^0$ continuous-time linear cocycle in the sense that $\w\mapsto\Phi_A(t,\w)$ is continuous for each $t\in\R$. 
\end{example}

\medskip

\begin{example}\label{Linfty} \textbf{Bounded linear differential systems:} Let $L^\infty(M,\mathbb{R}^{2\times 2})$ be the space of (essentially bounded) measurable matrix-valued maps $A\colon M\rightarrow \mathbb{R}^{2\times 2}$, satisfying
$$\|A\|_\infty:=\operatorname{ess}\sup_{\w\in M} \|A(\w)\|<\infty,$$
where  $\|\cdot\|$ denotes de standard Euclidean matrix norm. Since $M$ is compact the volume measure $\mu$ is finite and so $A\in L^1(\mu)$. It follows from \cite[Thm. 2.2.8]{A} (see also Example 2.2.8 in this reference) that if $A\in L^1(\mu)$ then ~\eqref{LVE10} generates a unique (up to indistinguishability) weak solution $\Phi_A$. Moreover, for each $t\in\R$, $\w\mapsto\Phi_A(t,\w)\in L^\infty(M,\R^{2\times2})$; see \cite{CongSon}.

\end{example}

Now we give a motivation to see equation \eqref{SOLH} as a linear differential system that can be studied from the asymptotic point of view. Let $\mathpzc{K}\subset \mathbb{R}^{2\times 2}$ be the set of matrices of type
	$$
		\begin{pmatrix} 0&1\\ b & a \end{pmatrix}
$$
	for real numbers $a, b$. This set $\mathpzc{K}$ define an affine space of $\mathbb{R}^{2\times 2}$ in a sense that 
	$$\mathpzc{K}=\left\{ K\in \mathbb{R}^{2\times 2}\colon K=U+V\right\},$$
	where $U=\left(\begin{matrix} 0&1\\ 0 &0 \end{matrix}\right)$ and $V=\left(\begin{matrix} 0&0\\ x_v &y_v \end{matrix}\right)$ belongs to the 2-dimensional subvector space of $\mathbb{R}^{2\times 2}$ defined by 
	\begin{equation}\label{00}
	\mathpzc{V}=\left\{ V\in \mathbb{R}^{2\times 2}\colon V=\left(\begin{matrix} 0&0\\ x &y \end{matrix}\right)\,\,\,\text{where}\,\,\,x,y\in \mathbb{R}\right\},
	\end{equation}
	and so $\mathpzc{K}$ is a $\mathpzc{V}$-torsor.

	Consider measurable and $L^\infty$ ($C^0$ will also be considered) maps $\alpha\colon M \to \R$ and $\beta\colon M \to \R$ and the second order linear differential equation based on the flow $\varphi^t$ given by
	\begin{equation}\label{dampedp}
		\ddot{x}(t) + \alpha(\varphi^{t}(\omega)) \dot{x}(t) + \beta(\varphi^{t}(\omega)) x(t) = 0.
	\end{equation}
		Taking $y(t)= \dot{x}(t)$ we may rewrite \eqref{dampedp} as the following vectorial linear system
	\begin{eqnarray}\label{E1}
		\dot{X}= A(\varphi^{t}(\omega))\cdot X,
	\end{eqnarray}
	where $X=X(t)=(x(t),y(t))^\top=(x(t),\dot x(t))^\top$ and $A\in L^\infty(M,\mathpzc{K})\subset L^\infty(M,\mathbb{R}^{2\times 2})$ is given by
	
	\begin{equation}\label{damp2}
\begin{array}{cccc}
A\colon &M & \longrightarrow & \mathbb{R}^{2\times2} \\& \w & \longmapsto &  \left(\begin{matrix}0&1\\ -\beta(\omega)& -\alpha(\omega)\end{matrix}\right)
\end{array}
\end{equation}
It follows from what we saw in Example~\ref{Linfty} that \eqref{LVE10} generates a Carath\'eodory solution $\Phi_A$. Given the initial condition $X(0)=v=(x(0),\dot x(0))$, the solution of \eqref{E1} is $X(t)=\Phi_A(t,\w)v$.

	We call the set $L^\infty(M,\mathpzc{K})\subset L^\infty(M,\mathbb{R}^{2\times 2})$ the \emph{$L^\infty$ kinetic cocycles} and $C^0(M,\mathpzc{K})\subset L^\infty(M,\mathpzc{K})$ the \emph{$C^0$ kinetic cocycles}. Two subsets of $L^\infty(M,\mathpzc{K})$
 will be of future interest: 
	\begin{itemize}
	\item the \emph{traceless/frictionless} ones characterized by $\alpha=0$. The $L^\infty$ or the $C^0$ kinetic cocycles will be  denoted by $L^\infty(M,\mathpzc{K}_{\,\,\Delta=0})$ or $C^0(M,\mathpzc{K}_{\,\,\Delta=0})$, respectively, and 
	\item the \emph{dissipative} ones characterized by $\alpha<0$. The $L^\infty$ or the $C^0$ kinetic cocycles will be denoted by  $L^\infty(M,\mathpzc{K}_{\,\,\Delta<0})$ or $C^0(M,\mathpzc{K}_{\,\,\Delta<0})$, respectively.
	\end{itemize}			

As we already discuss $\mathpzc{K}$ is not a vector subspace, but an affine subspace. In the sequel we will consider such perturbations $H\in \mathpzc{V}$ defined in \eqref{00} which we denote by $\mathbb{R}^{0\times 2}$. In conclusion, and since $\mathpzc{K}$ is a $\mathpzc{V}$-torsor, the set $\mathpzc{K}$ is closed under the sum of elements in the additive group $(\mathbb{R}^{0\times 2},+)$. In other words $(\mathbb{R}^{0\times 2},+)$ acts transitively say given any $K_1,K_2\in \mathpzc{K}$, there exists a unique $V\in \mathbb{R}^{0\times 2}$ such that $K_1+V =K_2$.

\subsection{Topologization of linear kinetic cocycles}

 We endow
$L^{\infty}(M,\mathbb{R}^{2\times 2})$ with the metric defined by
\begin{equation*}\label{Lin}
\rho_\infty(A,B)=\|A-B\|_{\infty},
\end{equation*}
where $A,B\in L^{\infty}(M,\mathbb{R}^{2\times 2})$. We also endow
$C^{0}(M,\mathbb{R}^{2\times 2})$ with the metric defined by
\begin{equation*}\label{C0n}
\rho_0(A,B)=\|A-B\|_{0},
\end{equation*}
where $A,B\in C^{0}(M,\mathbb{R}^{2\times 2})$ and $\|A\|_0={\underset{\w\in{M}}{\max}}\|A(\w)\|$.

\begin{proposition}\label{complete}
\noindent
\begin{enumerate}
\item[(i)] $(L^{\infty}(M,\mathbb{R}^{2\times 2}),\rho_\infty)$ and  $(C^{0}(M,\mathbb{R}^{2\times 2}),\rho_0)$ are complete metric spaces and, therefore, Baire spaces;
\item[(ii)] 
 $L^{\infty}(M,\mathpzc{K})$, $L^\infty(M,\mathpzc{K}_{\,\,\Delta=0})$ and $L^\infty(M,\mathpzc{K}_{\,\,\Delta<0})$ are $\rho_\infty$-closed;
 \item[(iii)]$C^{0}(M,\mathpzc{K})$, $C^0(M,\mathpzc{K}_{\,\,\Delta=0})$ and $C^0(M,\mathpzc{K}_{\,\,\Delta<0})$ are $\rho_0$-closed.
\end{enumerate}

\end{proposition}	

\begin{proof}
(i) We consider this property on each entry of the matrix. As $M$ is $\sigma$-finite we get that $L^{\infty}$ is the dual of $L^1$. Now we use the well-known result which states that the dual of any normed space is complete. The Baire property follows from Baire's category theorem. The set of bounded functions on $M$ endowed with the $C^0$ norm is complete. Moreover, the space $C^0$ functions on $M$ is a closed subset of the space bounded functions on $M$ and so is complete. (ii) and (iii) are trivial by the definition of $\mathpzc{K}$, $\mathpzc{K}_{\,\,\Delta=0}$ and $\mathpzc{K}_{\,\,\Delta<0})$. 
\end{proof}

\bigskip

\begin{remark} As we will see Theorem~\ref{tps4} is the $L^\infty$ version of Theorem \ref{tps}.  Clearly, Theorems \ref{tps2}, \ref{tps3} and \ref{tps5} also hold if we replace $C^0$ by $L^\infty$ in the statements \emph{mutatis mutandis}. \end{remark}

\bigskip
				
\subsection{Lyapunov exponents}
Take $A\in \mathbb{R}^{2\times 2}$, a vector $v\in\mathbb{R}^2$ and the solution $\Phi^t_A$ of $\dot U(t)=A\cdot U(t)$. The asymptotic growth of the expression $\frac{1}{t}\log\|\Phi^t_A\cdot v\|$ turns out to be a simple exercise in linear algebra. Simple calculations allow us to determine the spectral properties defined by eigenvalues and eigendirections. Fixing the terminology for what follows the logarithm of eigenvalues are called \emph{Lyapunov exponents} and the eigendirections are called \emph{Oseledets directions} which we will see in detail in a moment. Another very interesting problem is determining the stability of $\dot U(t)=A\cdot U(t)$ which aims to establish whether the qualitative behaviour remains similar when we perturb $A$. This issue is well understood in the autonomous case. In particular if $A\in \mathpzc{K}$ the asymptotic stability analysis of second order homogeneous differential equation with constant parameters $\alpha$ and $\beta$ is a problem already studied and understood. A considerably more complicated situation worthy of in-depth study was considered in Lyapunov's pioneering works and aimed at considering the non-autonomous case $\dot U(t)=A(t)\cdot U(t)$, where $A$ is a matrix depending continuously on $t$.  Not only the asymptotic demeanor of $\frac{1}{t}\log\|\Phi^t_A\|$ as well as its stability reveals itself as a substantially harder question. At this point the reader must agree that when $A(t)\in \mathpzc{K}$ for all $t$ the problem is associated to the asymptotic stability analysis of second order homogeneous differential equation with varying parameters $\alpha(t)$ and $\beta(t)$.

Given $A\in L^{\infty}(M,\mathbb{R}^{2\times 2})$ then Oseledets theorem (see e.g. \cite{O,A, JPS}) guarantees that for $\mu$ almost every $\omega\in M$, there exists a $\Phi_A(t,\w)$-invariant splitting called \emph{Oseledets splitting} of the fiber $\mathbb{R}^{2}_{\omega}=E^{1}_{\omega}\oplus E^{2}_{\omega}$ and real numbers called \emph{Lyapunov exponents} $\lambda_{i}(A,\omega)$, $i=1,2$, such that:
	
	\begin{equation*}\label{limit}
\lambda_{i}(A,\omega)=\lambda(A,\w,v^i)=	\underset{t\rightarrow{\pm{\infty}}}{\lim}\frac{1}{t}\log{\|\Phi(t,\omega) v^{i}\|}
	\end{equation*}
	for any $v^{i}\in{E^{i}_{\w}\setminus\{\vec{0}\}}$ and $i=1, 2$. If we do not count the multiplicities, then we have $\lambda_{1}(A,\omega)\geq \lambda_{2}(A,\omega)$. Moreover, given any of these subspaces $E^{1}_{\w}$ and $E^{2}_{\w}$, the angle between them along the orbit has subexponential growth, meaning that
	\begin{equation*}
	\lim_{t\rightarrow{\pm{\infty}}}\frac{1}{t}\log\sin(\measuredangle(E^{1}_{\varphi^{t}(\omega)},E^{2}_{\varphi^{t}(\omega)}))=0.	
\end{equation*}
		If the flow $\varphi^{t}$ is ergodic, then the Lyapunov exponents and the dimensions of the associated subbundles are constant $\mu$ almost everywhere and we refer to the former as $\lambda_1(A)$ and $\lambda_2(A)$, with $\lambda_1(A)\geq\lambda_2(A)$. We say that $A$ has \emph{simple} (Lyapunov) spectrum (respectively \emph{trivial} (Lyapunov) spectrum) if for $\mu$ a.e. $\w\in M$, $\lambda_{1}(A,\omega)> \lambda_{2}(A,\omega)$ (respectively $\lambda_{1}(A,\omega)=\lambda_{2}(A,\omega)$). For details on these results on linear differential systems see \cite{A} (in particular, Example 3.4.15). 
				
\medskip

\section{Continuous dependence and perturbations}

We define the \emph{integrated top Lyapunov exponent function} of the 
system $A$, over any measurable, $\varphi^{t}$-invariant
set $\Gamma\subseteq{M}$ by:
\begin{equation}\label{uscf}
\begin{array}{cccc}
LE(\cdot,\Gamma)\colon& L^{\infty}(M,\mathbb{R}^{2\times 2}) & \longrightarrow & \mathbb{R} \\
& A & \longmapsto & \displaystyle\int_{\Gamma}\lambda^{+}(\omega)\,d\mu(\omega)
\end{array}
\end{equation}
The main goal now is to understand if $LE(\cdot,\Gamma)$ display nice continuity properties. For that purpose next simple result gives continuous dependence of the solution on the infinitesimal generator and will be useful to prove Lemma~\ref{usc} below when we get that the function $LE(\cdot,\Gamma)$ is upper semicontinuous.

\begin{lemma}\cite{CongSon}\label{lemma: hat A}
Let $A\in L^\infty(M,\mathbb R^{2\times2})$, $\alpha,\beta\in\R$, $\alpha\leq\beta$, and $\hat A >\|A\|_\infty$. There exists a $\mu$-full subset $\tilde M\subset M$ (depending on $\hat A$) $\varphi^t$-invariant ($\varphi^t(\tilde M)\subset \tilde M$, $\forall t\in\R$), such that for all $\omega\in \tilde M$ we have

\begin{equation*}\label{eq:hatA}
\int_\alpha^\beta \|A(\varphi^s(\w))\|\,ds\leq\hat A(\beta-\alpha).
\end{equation*}
\end{lemma}
By Gronwall's inequality we have the following.
\begin{corollary}
\label{cor: Phi^t leq exp hat A t}
\,
\begin{enumerate}
\item Let $A\in C^0(M,\mathbb R^{2\times2})$. For almost every $\omega\in M$ and  every $t\in\R$ one has
\[
\|\Phi_A(t,\w)\| \leq \exp(|t|\,\|A\|_0).
\]
\item Let $A\in L^\infty(M,\mathbb R^{2\times2})$. For almost every $\omega\in M$,  every $t\in\R$ and every $\hat A>\|A\|_\infty $ one has
\[
\|\Phi_A(t,\w)\| \leq \exp(|t|\,\hat A)
\]
\end{enumerate}
In particular, if  $A\in C^0(M,\mathbb R^{2\times2})$ or $A\in L^\infty(M,\mathbb R^{2\times2})$ then, for almost every $\w$ and every fixed $t\in\R$ we have $\Phi_A^t= \Phi_A(t,\cdot)\in C^0(M,\mathbb R^{2\times2})$ or $\Phi_A^t\in L^\infty(M,\mathbb R^{2\times2})$, respectively.
\end{corollary}

\begin{lemma}\label{nota1}
Let $A\in L^\infty(M,\mathbb{R}^{2\times 2})$.  There exists a $\mu$-full subset $\tilde M\subset M$,  $\varphi^t$-invariant, such that for  any $\epsilon>0$, there exists $\tilde\tau>0$ such that for all $\omega\in \tilde M$  and  $s\in[0,\tilde\tau[$,
writting 
$$\Phi_A(t,\w) =\begin{pmatrix}\varphi_{11}(t,\w)&\varphi_{12}(t,\w)\\ \varphi_{21}(t,\w) &\varphi_{22}(t,\w) \end{pmatrix}$$
 we have
\begin{equation*}\label{eq:Phi_Id}
|\varphi_{ii}(s,\w)-1|<\epsilon\,\,\,\text { and } \,\,\,|\varphi_{ij}(s,\w)|<\epsilon,
\end{equation*}
for all $i,j=1,2$, $i\neq j$, and $0\leq s \leq\hat\tau$. 
\end{lemma}

\begin{proof}
Fix any $\hat A>\|A\|_\infty$ and consider $\tilde M\subset M$ given by Lemma~\ref{lemma: hat A}.  Let $\hat\tau$ be such that ${\hat{\tau}}^2\hat A  e^{\hat A\hat\tau}<\epsilon$. From~\eqref{eq:LDS0}, Lemma~\ref{lemma: hat A} and Corollary~\ref{cor: Phi^t leq exp hat A t} we have for all $\w\in \tilde M$
\begin{eqnarray*}
\left\|\Phi_A(\hat\tau,\w)-\text{Id}\,\right\|
&\leq &
\int_0^{\hat\tau}\|A(\varphi^s(\w))\|\cdot\|\Phi_A(s,\w)\|\,ds\leq 
\int_0^{\hat\tau}\hat A \hat\tau e^{\hat A\hat\tau}\,ds
\leq \epsilon.
\end{eqnarray*}
\end{proof}

\bigskip

Now we are in position to obtain the continuous dependence.

\begin{lemma}\label{contidep}
Let $A,B\in L^{\infty}(M,\mathbb{R}^{2\times 2})$. For any fixed $t>0$ we have
\begin{equation}\label{liminf}
\underset{\rho_\infty(A,B)\rightarrow 0}{\lim}\,\rho_\infty(\Phi^t_A,\Phi^t_B)=0.
\end{equation}
\end{lemma}
\begin{proof}
Since $(A-B)\in L^{\infty}(M,\mathbb{R}^{2\times 2})$, from Lemma~\ref{lemma: hat A}, for any $\epsilon>0$ we have
\[
\int_0^t \|A(\varphi^s(\w))-B(\varphi^s(\w))\|\,ds\leq (\rho_\infty(A,B)+\epsilon)t.
\]
Take $\hat A>\|A\|_\infty$ and $\hat B>\|B\|_\infty$. 
From the Carath\'eodory solution \eqref{eq:LDS0} and fixing $\w\in M$,  taking $u(t)=\|\Phi_{A}(t,\omega)-\Phi_{B}(t,\omega)\|$,  we have:
\begin{eqnarray*}
u(t)&=&\left\|\int\limits_{0}^{t}A(\varphi^{s}(\omega))\Phi_{A}(s,\omega)-B(\varphi^{s}(\omega))\Phi_{B}(s,\omega)\,ds\right\|\\
&\leq&\int\limits_{0}^{t}\left\|A(\varphi^{s}(\omega))(\Phi_{A}(s,\omega)-\Phi_{B}(s,\omega))\right\|\,ds+\int\limits_{0}^{t}\left\| (A(\varphi^{s}(\omega))-B(\varphi^{s}(\omega)))\Phi_{B}(s,\omega)\right\|\,ds\\
&\leq&\int\limits_{0}^{t}\underbrace{\|A(\varphi^{s}(\omega)\|}_{\beta(s)}\,\underbrace{\|\Phi_{A}(s,\omega)-\Phi_{B}(s,\omega)\|}_{u(s)} \,ds+\underbrace{(\rho_\infty(A,B)+\epsilon)t\,\text{exp}(t\hat B)}_{\alpha(t)}.
\end{eqnarray*}

In summary, we have
$$u(t)\leq \int_0^t \beta(s)u(s)\,ds+\alpha(t).$$
Using Gr\"onwall's inequality we get:
$$u(t)\leq \alpha(t)\exp \left(\int_0^t \beta(s)\,ds\right),$$
that is
\begin{align*}
\|\Phi_{A}(t,\omega)-\Phi_{B}(t,\omega)\|
&\leq 
(\rho_\infty(A,B)+\epsilon)t\,\text{exp}(t\hat B)\text{exp}\left(\int\limits_{0}^{t}\|A(\varphi^{s}(\omega))\|\,ds\right)\\
&\leq 
(\rho_\infty(A,B))+\epsilon)t\,\text{exp}(t\hat B)\text{exp}(t\hat A)
\end{align*}
and \eqref{liminf} holds since $\epsilon>0$ is arbitrary.
\end{proof}

\begin{lemma}\label{usc}
The function $LE(\cdot,\Gamma)$ is upper semicontinuous.
\end{lemma}

\begin{proof}
Take $A\in \mathpzc{K}$ and let $a_n(A)=\int_{\Gamma}\log\|\Phi_A(n,\omega)\|\,d\mu(\omega)$. Clearly, the sequence $a_n$ is subadditive. Moreover, by Fekete's lemma
$$LE(\cdot,\Gamma)=\underset{n\rightarrow+\infty}{\lim}\frac{a_n(A)}{n}=\underset{n\in\mathbb{N}}{\inf}\frac{a_n(A)}{n}.$$
Now Lemma~\ref{contidep} and the continuity of the norm and of the logarithm allow us to obtain that the map $A\mapsto a_n(A)$ is continuous. Using the subadditivity of the norm we obtain:
$$LE(A,\Gamma)=\underset{n\in\mathbb{N}}{\inf}\,\frac{1}{n}\int_{\Gamma}\text{log}\|\Phi_A(n,\w)\|d\mu(\w).$$
Since $LE(\cdot,\Gamma)$ is the infimum of a sequence of continuous functions it is upper semicontinuous. 

\end{proof}

\medskip

\section{The toolbox of kinetic perturbations} \label{toolbox}

\subsection{Accessibility and kinetic perturbations} 

We now provide some perturbative tools that allow us to rotate directions defined by Oseledets fibres and within the class of linear (conservative) kinetic cocycles. We'll take care to make the perturbations conservative so that we can use them in more restricted contexts. These perturbations will be key to proving our results. In \cite[Theorem 5]{BV2} was proved that discrete accessible families of cocycles satisfy the same conclusions of the theorems of the present paper. So, our main goal is to prove that families of kinetic cocycles are accessible. Accordingly, next definition is central in our study:

\begin{definition}\label{acc} (\textbf{Accessibility})
We say that the set $\mathcal{G}\subseteq L^\infty(M,\mathbb{R}^{2\times 2})$ is \emph{accessible} if for all $\epsilon>0$, there exists $\theta>0$ with the following properties: Given $A\in \mathcal{G}$, a $\mu$-generic point $\w\in M$ and $u,v$ in the projective space $\mathbb{R}P^1_\w$ with\footnote{We consider $\measuredangle(u,v)$ as the minimum angle between the half-lines generated by $u$ and $v$.} $\measuredangle(u,v)<\theta$, then there exists $B\in \mathcal{G}$ with $\|A-B\|_\infty<\epsilon$, such that $\Phi_B(1,\w) u=\Phi_A(1,\w) v$. 
\end{definition}

\begin{remark}
Definition~\ref{acc} implies that when considering the full $\mu$-measure set $$M_A=\{\w\in M\colon \|A(\w)\|\leq \|A\|_\infty\}$$ the set $$\{\Phi_A(1,\w)\colon A\in \mathcal{G},\w\in M_A\}\subset \text{GL}(2,\mathbb{R})$$ is accessible in the sense of \cite[Definition 1.2]{BV2} with $\nu=1$ and $C=C(\|A\|_\infty)$.
\end{remark}

Definition~\ref{acc} deals with a perturbation $B$ of a cocycle $A$. Now we define what we mean by perturbation. Take $A\in L^{\infty}(M,\mathpzc{K})$ and a non-periodic point $\omega\in M$. We will define a perturbation $B\in L^{\infty}(M,\mathpzc{K})$ of $A$ supported on a time-$\tau$ segment of orbit, $\tau>0$, starting on $\varphi^T(\w)$, for some $T\geq0$ and $\w\in M$:
$$\mathpzc{S}=\varphi^{[T,T+\tau]}(\w)=\bigl\{\varphi^{s}(\varphi^T(\w)): s\in[0,\tau]\bigr\}.$$
\begin{definition}\label{Defpert} (\textbf{Perturbations}) Take $B_0\in L^{\infty}(\mathcal{S},\mathpzc{K})$. The \emph{perturbation $B\in L^{\infty}(M,\mathpzc{K})$ of a given $A\in L^{\infty}(M,\mathpzc{K})$ supported on 
$\mathpzc{S}$} is defined  by:
		\begin{equation}\label{flow}B(\w):=\left\{ \begin{array}{lll}
	A(\w)& \mbox{if $\w\notin \mathpzc{S}$}\\ B_0(\w)& \mbox{if $\w\in \mathpzc{S}$}\end{array}\right.
\end{equation}
\end{definition}

Observe that even when $A\in C^{0}(M,\mathpzc{K})$ and $B_0\in C^{0}(\mathcal{S},\mathpzc{K})$, $B$ is not necessarily continuous. Nevertheless, $t\mapsto \Phi_B(t,\w)$ is continuous.

\medskip

\subsection{Why should we rotate solutions?}
The goal of perturbations as in Definition~\ref{acc} is to to cause some rotational effect. The strategy to obtain trivial Lyapunov spectrum under a small perturbation relies on a fundamental idea which goes back to the works of Novikov \cite{Nov} and Ma\~n\'e \cite{M1}: \emph{Rotate Oseledets directions} an idea which nowdays as synonyms in the literature like \emph{twisting} or \emph{accessibility}. The main idea is quite intuitive: For $i=1,2$, the Lyapunov exponent along $E_\w^i$ represents an assymptotic average $\lambda_i(A,\w)$. It seems plausible that combining these directions, say rotate one into another, cause a mix of the two rates. For example, if for a huge iterate $\tau$ we spend $\tau/2$ on a rate 
$$\left\|\Phi_A({\tau}/{2},\w)\Big|_{E^1_\w}\right\|\sim e^{\frac{\tau}{2}\lambda_1}$$ and the remaining $\tau/2$ on a rate 
$$\left\|\Phi_A({{\tau}/{2}}, \varphi^{{\tau}/{2}}(\w))\Big|_{E^2_{ \varphi^{{\tau}/{2}}(\w)}}\right\|\sim  e^{\frac{\tau}{2}\lambda_2}.$$
At the end of the day we get a rate 
$$\left\|\Phi_B({\tau},\w)\right\|\sim  e^{\tau\left(\frac{\lambda_1+\lambda_2}{2}\right)},$$
for the perturbed cocycle $B$. But we can't just suddenly change directions like that, right? In fact we will consider somewhere in the middle of the trajectory a segment of orbit of size 1 and on that segment we will swap the directions. This approach was successfully put into practice in the works \cite{B,BV2, Be} considering $C^0$ cocycles  on the special linear group, in the $L^p$ special linear group case by \cite{BVi}, in the $L^p$ kinetic case in \cite{ABV} and were inspired on the ideas contained in \cite{M1,M2}. However, we need to make sure that we are able to do this in the world of kinetic systems. This class has very particular idiosyncrasies (cf. \cite{ABV,ABV2, Be4}). So we begin the carefully construct of the perturbations.

\subsection{Escaping vertical directions}\label{esc}

Given $A\in L^{\infty}(M,\mathpzc{K})$ the perturbation $B\in L^{\infty}(M,\mathpzc{K})$ of $A$ as \eqref{flow} in Definition~\ref{Defpert} will be defined by $B=A+H$ where $H\in L^{\infty}(M,\mathbb{R}^{0\times 2})$. So we begin by considering a constant infinitesimal generator $H\in L^{\infty}(M,\mathbb{R}^{0\times 2})$ given by
 \begin{equation}\label{H00}
 H(\w)=\begin{pmatrix}0&0\\ -\beta & -\alpha \end{pmatrix},
 \end{equation}
with $\alpha,\beta\in\R$. Clearly, $H$ acts in $u=(u_1,u_2)\in \mathbb{R}^2$ as the vector field 
 \begin{equation*}\label{oH}
 H\cdot u=(0,-\beta u_1-\alpha u_2),
 \end{equation*}
 say just like a vertical vector field. The next observation shows that vectors that are nearly vertical are difficult to rotate to the vertical direction by the action of $H\in \mathbb{R}^{0\times 2}$. This type of behavior was already observed in \cite[Example 5]{BV2} when related to Schr\"odinger cocycles which are a particular subclass of discrete $\text{SL}(2,\mathbb{R})$ cocycles associated to our kinetic ones. Indeed, Schr\"odinger cocycles arises from a $\mathpzc{V}$-torsor defined as in \eqref{00} for a particular choice of
 $$U=\left(\begin{matrix} 0&0\\ 1 &0 \end{matrix}\right)\,\,\,\text{ and }\,\,\,V=\left(\begin{matrix}E-Q(\w) &-1\\ 0&0\end{matrix}\right)$$ where $E$ is an energy and $Q$ a potential cf. \cite[Example 2]{BV2}). See also \S\ref{Sc} where a continuous-time counterpart is discussed.

\begin{remark}
Let $\epsilon>0$ and $u=(\cos\theta,\sin\theta)$ where $\theta\in \left]\frac{\pi}{2}-\epsilon,\frac{\pi}{2}\right[$. Let also $H$ be as in \eqref{H00} with $\alpha,\beta<0$. Let $H\cdot u=(0,\kappa)$, for $\kappa=-\beta u_1-\alpha u_2>0$. The angle $\gamma=\gamma(\theta,\kappa)$ between $u$ and $u+H\cdot u$ is such that:
$$\gamma=\arccos\frac{(\cos\theta,\sin\theta)\cdot(\cos\theta,\sin\theta+\kappa)}{\|(\cos\theta,\sin\theta+\kappa)\|}=\arccos\frac{1+\kappa\sin\theta}{\sqrt{1+2\kappa\sin\theta+\kappa^2}}\underset{\theta\rightarrow\frac{\pi}{2}}{\longrightarrow}\arccos1=0.$$
\end{remark}

Fortunately, as we will prove in Lemma~\ref{EsLem}, vectors tend to steer away from the vertical position as the cocycles in our consideration take effect. In other words our kinetic cocycles are transversal to the vertical direction. Indeed, the infinitesimal generator $A\in L^{\infty}(M,\mathpzc{K})$, with
 \begin{equation*}\label{A00}
 A(\w)=\begin{pmatrix}0&1\\ -\beta(\w) & -\alpha(\w) \end{pmatrix}
 \end{equation*}
 acts in $u=(0,1)\in \mathbb{R}^2$ as the vector field 
 \begin{equation*}\label{oA}
A(\w)\cdot u=(1,-\alpha(\w)),
 \end{equation*}
and since $A$ is bounded, then $\alpha$ is bounded and so $(1,-\alpha)$ is transversal to $(0,1)$.

We introduce now some useful notation on cones. Given $\gamma>0$ let $\mathcal{C}_\gamma\subset\mathbb{R}^2$ denote the \emph{vertical cone} defined by 
$$(v_1,v_2)\in \mathcal{C}_\gamma\,\,\,\text{if}\,\,\, |v_1|<\gamma|v_2|.$$

Next result says that kinetic cocycles tend to escape from the vertical direction and no perturbation is needed. The effect of a kinetic infinitesimal generator is to escape the cone $\mathcal{C}_\gamma$. At $\nu\in \mathcal{C}_\gamma$ the vector field pushes in the direction $A\cdot\nu$ helping to get $\nu$ out of the cone.

\begin{figure}

\begin{center}
\begin{tikzpicture}

    \draw[gray] (-0.5,-4) -- (0.5,4);
    \draw[gray] (-0.5,4) -- (0.5,-4);

    \fill[gray!10] (-0.5,-4) -- (0.5,4) -- (-0.5,4) -- (0.5,-4) -- cycle;

    \node[gray, below right] at (-0.6,4.1) {$\mathcal{C}_\gamma$};
     \node[gray, below right] at (0.9,4.2) {$A\cdot\nu$};
      \node[gray, below right] at (-0.1,3.4) {$\nu$};

        \draw[->] (-4,0) -- (4,0) node[right] {$x$};
    \draw[->] (0,-4) -- (0,4.2) node[above] {$\dot x$};
    
     \draw[->, thick] (0,0) -- (0.2,3);
          \draw[->, thick] (0.2,3) -- (1,4);
     \draw[->, thick, dashed] (0.2,3) --  (0.8,2.82);

\end{tikzpicture}
\end{center}
\caption{Escaping from thin vertical cones.}\label{Fig1}
\end{figure}
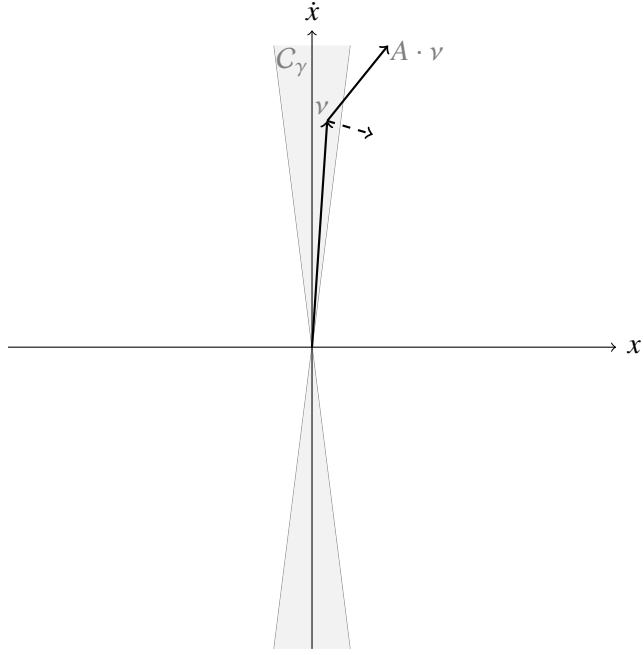

\begin{lemma}\label{EsLem} \textbf{Escaping lemma:}
	Let be given $A\in L^{\infty}(M,\mathpzc{K})$. There exist $\hat\tau\in\left]0,\frac{1}{2}\right[$ and  $\gamma\in]0,1[$ such that for almost every point $\w\in M$ and $v\in \mathcal{C}_\gamma\subset\mathbb{R}^2_\w$ we have $\Phi_{A}(\hat\tau,\omega)v\notin \mathcal{C}_\gamma\subset\mathbb{R}^2_{\varphi^{\hat\tau}(\w)}$. 
\end{lemma}

\begin{proof}
Consider $A\in L^{\infty}(M,\mathpzc{K})$ given as
\begin{equation*}\label{A00}
 A(\w)=\begin{pmatrix}0&1\\ -\beta(\w) & -\alpha(\w) \end{pmatrix}
 \end{equation*}
and fix 
\begin{equation}\label{min0}
\ell>\max\{\text{ess\,sup} \,\alpha,\text{ess\,sup} \,\beta\}.
\end{equation}
For simplicity, let us write
\[
\Phi_A(t,\w) =\begin{pmatrix}\varphi_{11}^t&\varphi_{12}^t\\ \varphi_{21}^t &\varphi_{22}^t \end{pmatrix}.
\]
Consider any $0<\epsilon<\frac12$ and let $\hat{\tau}$ be such that for almost every $\w$,
\begin{equation}\label{min1}
|\varphi_{ii}^s-1|<\epsilon\,\,\,\text { and } \,\,\,|\varphi_{ij}^s|<\epsilon,
\end{equation}
for all $i,j=1,2$, $i\neq j$, and $0\leq s \leq\hat\tau$.
 Let $v\in \mathcal{C}_\gamma\subset\mathbb{R}^2_\w$. Without loss of generality we may assume that $v$ has the form $v=(\tilde \gamma,1)$, where $0<|\tilde\gamma|<\gamma$. 
 
 Now, when we apply the transformation $A(\w)$ to $v$, we obtain $$A(\w) \cdot v = (1, -\beta(\w) \tilde{\gamma} - \alpha(\w)).$$ 
We need to make sure that the infinitesimal action of $A$ along the orbit of $\omega$ over $v$ induces a horizontal translation on $v$ by escaping the vertical cone  $\mathcal C_\gamma$ (see Figure \ref{Fig1} above). To achieve this, it is necessary that $A(\w)\cdot v$ does not belong to $\mathcal{C}_\gamma$. The worst case scenario occurs when $\alpha=\beta=-\ell$ implying that $A(\w)\cdot v$ is $(1, \ell(\tilde{\gamma}+1))$. This means that the infinitesimal generator tends to force $\Phi_A(t,\w)$ move away $v$ from a suitable $\mathcal{C}_\gamma$ at leats for times near $t=0$. We need to guarantee that indeed we can choose a narrow cone such that all vectors inside leave the cone in a small fraction of time, making use of the effort in the horizontal given by the infinitesimal generator. 

Let us choose $\gamma>0$ small enough such that for all $0<\tilde\gamma<\gamma$
\begin{equation}\label{min}
\min\{|\tilde\gamma\pm\hat\tau\tilde\gamma\epsilon+(1\pm\epsilon)\hat\tau|\}>\gamma,
\end{equation}
where the minimum is taken over all the possible combinations of signs.
We have

\begin{eqnarray*}
\Phi_{A}(\hat\tau,\omega)v&=&v+\int_0^{\hat\tau} A(\varphi^t(\w))\Phi_{A}(t,\omega) v\,dt\\
&=&
\begin{pmatrix} \tilde\gamma\\ 1\end{pmatrix}+\int_0^{\hat\tau}  \begin{pmatrix} 0 & 1\\ -\beta(\varphi^t(\w)) & -\alpha(\varphi^t(\w))\end{pmatrix}
\begin{pmatrix}\varphi_{11}^s&\varphi_{12}^s\\ \varphi_{21}^s &\varphi_{22}^s \end{pmatrix}
\begin{pmatrix} \tilde\gamma\\ 1\end{pmatrix}\,dt\\
&=&
\begin{pmatrix} \tilde\gamma + \int_0^{\hat\tau} \tilde\gamma\varphi_{21}^s + \varphi_{22}^s\,ds \\ 
1 - \int_0^{\hat\tau} \beta(\varphi^t(\omega))(\tilde\gamma\varphi_{11}^s + \varphi_{12}^s) + \alpha(\tilde\gamma\varphi^t(\omega))\varphi_{21}^s+\varphi_{22}^s)\,ds \end{pmatrix}
=\begin{pmatrix} z_1\\ z_2\end{pmatrix}
\end{eqnarray*}
Now, from \eqref{min0}, \eqref{min1} and \eqref{min} we get
\begin{eqnarray*}
\frac{|z_1|}{|z_2|}\geq\frac{\min\{|\hat\gamma\pm\hat\tau\tilde\gamma\epsilon+(1\pm\epsilon)\hat\tau|\}}{1 + \hat\tau\ell(2|\tilde\gamma|\epsilon+2\epsilon+|\tilde\gamma|+1+\epsilon)}>\gamma.
\end{eqnarray*}
Hence, $\Phi_{A}(\hat\tau,\omega)v\notin \mathcal{C}_\gamma\subset\mathbb{R}^2_{\varphi^{\hat\tau}(\w)}$.

\end{proof}

\bigskip

\begin{remark}\label{rem3}
Notice that if we consider $\theta_0>0$ very small, then from the proof of Lemma~\ref{EsLem} if $u$ is such that $\measuredangle(u,v)<\theta_0$, then we also have $\Phi_{A}(\hat\tau,\omega)u\notin \mathcal{C}_\gamma\subset\mathbb{R}^2_{\varphi^{\hat\tau}(\w)}$. 

\end{remark}

\bigskip

\subsection{Conservative kinetic rotations}

We start by introducing a simple example where we can observe the flavour of a kinetic rotation seen as an action on a fixed fiber.
		
\begin{figure}
\begin{center}
\begin{tikzpicture}
    \draw[->] (-2,0) -- (2,0) node[right] {$x$};
    \draw[->] (0,-2) -- (0,2) node[above] {$\dot x$};
    \draw (-1,2) node[above] {$\eta>0$};
    \draw[->] (-0,1.2) -- (0.2,1.2) node[right] {{\tiny saddle-node}};
        \draw[->] (-0,-1.2) -- (0.2,-1.2) node[right] {{\tiny saddle-node}};

    \draw (0,0) circle [radius=1];

\draw[->,ultra thick] (0,0.2) ++(0:0.6) arc (0:60:0.6);
 \draw[->,ultra thick] (0,-0.2) ++(-60:0.6) arc (-60:0:0.6);
 \draw[->,ultra thick] (0,0.2) ++(120:0.6) arc (120:180:0.6);
  \draw[->,ultra thick] (0,-0.2) ++(180:0.6) arc (180:240:0.6);

  \draw[->,ultra thick] (0,0.2) ++(0:0.6) arc (0:60:0.6);
 \draw[->,ultra thick] (0,-0.2) ++(-60:0.6) arc (-60:0:0.6);
 \draw[->,ultra thick] (0,0.2) ++(120:0.6) arc (120:180:0.6);
  \draw[->,ultra thick] (0,-0.2) ++(180:0.6) arc (180:240:0.6);

\fill (0,1) circle [radius=2pt];
\fill (0,-1) circle [radius=2pt];
\end{tikzpicture}
\begin{tikzpicture}

    \draw[->] (-2,0) -- (2,0) node[right] {$ x$};
    \draw[->] (0,-2) -- (0,2) node[above] {$\dot x$};
    \draw (-1,2) node[above] {$\eta<0$};
    \draw[->] (-0,1.2) -- (0.2,1.2) node[right] {{\tiny saddle-node}};
            \draw[->] (-0,-1.2) -- (0.2,-1.2) node[right] {{\tiny saddle-node}};

    \draw (0,0) circle [radius=1];

\draw[->,ultra thick] (0,0.2) ++(60:0.6) arc (60:0:0.6);
 \draw[->,ultra thick] (0,-0.2) ++(0:0.6) arc (0:-60:0.6);
 \draw[->,ultra thick] (0,0.2) ++(180:0.6) arc (180:120:0.6);
  \draw[->,ultra thick] (0,-0.2) ++(240:0.6) arc (240:180:0.6);

\fill (0,1) circle [radius=2pt];
\fill (0,-1) circle [radius=2pt];
\end{tikzpicture}

\end{center}
\caption{Illustration for Example \ref{SM}.}\label{Fig2}
\end{figure}
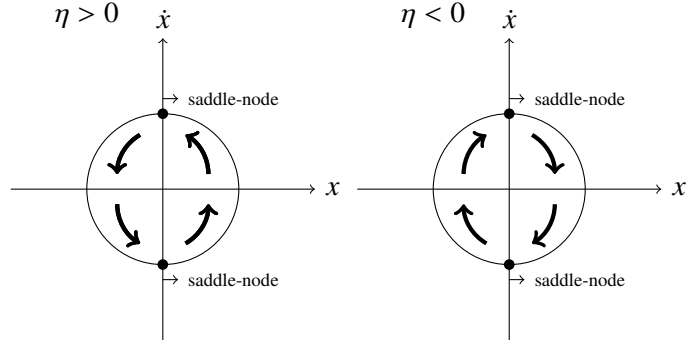

	\begin{example}\label{SM} \textbf{(shear matrix)}
	The fundamental solution $R^t\colon \mathbb{R}^2_\w\rightarrow \mathbb{R}^2_\w$ of the autonomous linear variational equation $\dot X(t)=\mathcal{A}\cdot X(t)$, where
	$$
\begin{array}{cccc}
\mathcal{A}\colon& \mathbb{R}^2_\w & \longrightarrow & \mathbb{R}^2_\w \\
& v & \longmapsto & \begin{pmatrix} 0 & 0\\ \eta & 0\end{pmatrix}v
\end{array}
$$
and $\eta\in\mathbb{R}$, is given by the shear matrix $$R^t= \begin{pmatrix} 1 & 0\\ \eta t & 1\end{pmatrix}$$ with shear parallel to the $y$ axis. Notice that from the linear variational equation we get $\dot R^t\cdot[R^t]^{-1}=\mathcal{A}$. The single eigenvalue of $R^t$ is $1$ with eigenvector $(0,1)$ which is a saddle-node fixed point in the projective space . If $\eta>0$ (respectively $\eta<0$) the dynamics in the projective space is a counter-clockwise `rotation' (respectively clockwise) except of course the direction fixed direction $(0,1)$ (see Figure~\ref{Fig2}). The angle between a vector and its image under $R^t$ decreases to 0 as we get close to the vertical direction inside a cone $C_\gamma$, and is bounded away from 0 outside the cone. Indeed, take $(\tilde\gamma,1)\in \mathcal{C}_\gamma$ where $\gamma>0$, that is $|\tilde\gamma|<\gamma$. Then, $R^t(\tilde\gamma,1)=(\tilde\gamma,\tilde\gamma\eta t+1)$. Clearly, 
$\measuredangle((\tilde\gamma,1),(\tilde\gamma,\tilde\gamma\eta t+1))\rightarrow 0$ when $\gamma\rightarrow0$ and $\measuredangle(u,R^tu)\geq \theta> 0$ for a fixed $\theta$ when $u\notin \mathcal{C}_\gamma$. 
\end{example}

We are interested in producing a perturbation of composition type.
\begin{equation}\label{compotype}
\Phi_B(\tau,\w)u=\Phi_A(\tau,\w)\cdot R^\tau u
\end{equation}
which is somehow related to the equality $\Phi_B(1,\w) u=\Phi_A(1,\w) v$ appearing in Definition~\ref{acc} because 
\begin{eqnarray*}
\Phi_{B}(1,\omega)u&=&\Phi_{B}(1-\tau,\varphi^{\tau}(\w))\Phi_{B}(\tau,\omega)u\overset{\eqref{compotype}}{=}\Phi_{B}(1-\tau,\varphi^{\tau}(\w))\Phi_A(\tau,\w)\cdot R^\tau u\\
&\overset{R^\tau u=v}{=}&\Phi_{A}(1-\tau,\varphi^{\tau}(\w))\Phi_A(\tau,\w)v=\Phi_{A}(1,\omega)v,
\end{eqnarray*}
So the next result (cf. Lemma~\ref{rot00}) takes care to consider simultaneously three key aspects for a perturbation: is kinetic, is frictionless and will allow a composition type perturbation.

\begin{definition}\label{rot00old}
	Let be given $A\in L^{\infty}(M,\mathpzc{K})$ and $\xi\in\mathbb{R}$. For all $\w\in M$ and $t\geq0$ we define
	\begin{equation*}\label{OA}
	\begin{array}{cccc}
\mathcal{A}(t)\colon& \mathbb{R}^2_\w & \longrightarrow & \mathbb{R}^2_\w \\
& v & \longmapsto & \begin{pmatrix} \nu(t) & \mu(t) \\ \eta(t)  & \zeta(t)  \end{pmatrix}\cdot v
\end{array}
\end{equation*}
satisfying
\begin{equation}\label{deeefH}
\mathcal{A}_\xi(t)=\mathcal{A}(t)= \begin{pmatrix} \nu(t) & \mu(t) \\ \eta(t)  & \zeta(t)  \end{pmatrix}=[\Phi_A(t,\w)]^{-1}\cdot \begin{pmatrix} 0 & 0\\ \xi & 0   \end{pmatrix}\cdot \Phi_A(t,\w).
\end{equation}
	\end{definition}

	\bigskip
	
	\begin{remark}\label{re:H}
	Notice that fixing $A\in L^{\infty}(M,\mathpzc{K})$, $\tau>0$ and a non-periodic $\w\in M$ and defining 
\begin{equation}\label{Hdef}
H(\varphi^t(\w))= \left\{
\begin{array}{ll}
     \begin{pmatrix} 0 & 0\\ \xi & 0   \end{pmatrix}= \Phi_A(t,\w)\cdot \mathcal{A}(t)\cdot [\Phi_A(t,\w)]^{-1} & \text{if}\,\,\,t\in[0,\tau] \\
     ~\\
      0 & \text{otherwise}
\end{array} 
\right. 
\end{equation}	
	 we have $A+H\in L^{\infty}(M,\mathpzc{K})$ and for
	$A\in L^{\infty}(M,\mathpzc{K}_{\,\,\Delta=0})$ we have $A+H\in L^{\infty}(M,\mathpzc{K}_{\,\,\Delta=0})$.
	\end{remark}

		\bigskip

	\begin{lemma}\label{rot00}
	Let be given $A\in L^{\infty}(M,\mathpzc{K})$ (respectively $A\in L^{\infty}(M,\mathpzc{K}_{\,\,\Delta=0})$), $\xi\in\R$, $\tau>0$ and a non-periodic $\w\in M$. Take $\mathcal{A}(t)$ as in Definition~\ref{rot00old}. Then $B=A+H$ as in \eqref{Hdef} belongs to $L^{\infty}(M,\mathpzc{K})$ (respectively $L^{\infty}(M,\mathpzc{K}_{\,\,\Delta=0})$). Moreover, for $t\in[0,\tau]$, we have 
\begin{equation}\label{composi}
\Phi_B(t,\w)=\Phi_A(t,\w)\cdot R^t
\end{equation}
where $R^t\colon \mathbb{R}^2_\w\rightarrow \mathbb{R}^2_\w$ is the solution of the linear variational equation $\dot X(t)=\mathcal{A}(t)\cdot X(t)$.

	\end{lemma}
	
	\begin{proof} 
Notice that $B$ is the perturbation of $A$ by $A+H$ supported on $\varphi^{[0,\tau]}(\w)$. That $B\in L^{\infty}(M,\mathpzc{K})$ (respectively $B\in L^{\infty}(M,\mathpzc{K}_{\,\,\Delta=0})$) follows from Definition~\ref{rot00old} and Remark~\ref{re:H}.
	
	Noticing that $\dot R^t=\mathcal{A}(t)\cdot R^t$ we get by integrating by parts

\begin{eqnarray*}
	\int_{0}^{t} B(\varphi^{s}(\omega))\cdot\Phi_A(s,\w)\cdot R^sds&\overset{\eqref{Hdef}}{=}&\int_{0}^{t} [A(\varphi^{s}(\omega))\\
	&+&\Phi_A(s,\w)\cdot \mathcal{A}(s)\cdot [\Phi_A(s,\w)]^{-1}]\cdot\Phi_A(s,\w)\cdot R^sds\\
	&=&\int_{0}^{t} A(\varphi^{s}(\omega))\cdot\Phi_A(s,\w)\cdot R^sds+\int_{0}^{t} \Phi_A(s,\w)\cdot \dot R^sds\\
	&=&\Phi_A(s,\w) R^s\Big|_0^t-\int_0^t \Phi_A(s,\w)\cdot\dot R^s ds +\int_{0}^{t} \Phi_A(s,\w)\cdot \dot R^sds\\
	&=&\Phi_A(t,\w)R^t-\text{Id}\\	
	\end{eqnarray*}	
and so
	$$\Phi_A(t,\w)\cdot R^t=\text{Id}+\int_{0}^{t} B(\varphi^{s}(\omega))\cdot (\Phi_A(s,\w)\cdot R^s)ds$$
	Moreover we also have
	$$\Phi_B(t,\w)=\text{Id}+\int_{0}^{t} B(\varphi^{s}(\omega))\cdot\Phi_{B}(s,\omega) ds$$
and so \eqref{composi} follows from uniqueness of solutions of a linear integral equation 
$$U(t,\w)=\text{Id}+\int_{0}^{t} B(\varphi^{s}(\omega))\cdot U(s,\w)\,ds$$
when fixed the initial condition, and the lemma is proved.
		
	\end{proof}
	
	\bigskip

\subsection{Proving accessibility}
We start by showing that, for $t$ small, the action $R^t$ on a fixed fiber given by the solution of $\dot X(t)=\mathcal{A}(t)\cdot X(t)$, provides a rotational behavior stronger than the action $S^t$ of the solution of $\dot X(t)=H\cdot X(t)$, for a suitable constant $H$.

\begin{lemma}\label{nota2}
Let $A\in L^\infty(M,\mathbb{R}^{2\times 2})$, $\epsilon>0$ and $\gamma\in ]0,1[$. There exist $\tilde\tau\in \left]0,\frac{1}{2}\right[$ and $\theta=\theta(\tilde\tau)>0$ such that   for $\w\in M$ and $u=(\gamma,1)$ we have
\begin{equation}\label{first}
\max_{t\in [0,\tilde\tau]}\|H-\mathcal{A}_{\epsilon}(t)\|<\epsilon,
\end{equation}
and
\begin{equation}\label{second}
\measuredangle(u,R^{\tilde\tau} u)>\measuredangle(u,S^{\tilde\tau} u)=\theta
\end{equation}
where $$H=\begin{pmatrix} 0& 0 \\ \frac{\epsilon}{2} & 0 \end{pmatrix},$$ $R^t\colon \mathbb{R}^2_\w\rightarrow \mathbb{R}^2_\w$ is the solution of the linear variational equation $\dot X(t)=\mathcal{A}(t)\cdot X(t)$ for $\mathcal{A}(t)$ defined as in \eqref{deeefH} and $S^t\colon \mathbb{R}^2_\w\rightarrow \mathbb{R}^2_\w$ is the solution of the linear variational equation $\dot X(t)=H\cdot X(t)$.
\end{lemma}

\begin{proof}
We begin by considering Lemma~\ref{nota1} and choose $\tilde\tau\in \left]0,\frac{1}{2}\right[$ such that \eqref{first} holds. Recall the definitions of $H$ and $\mathcal{A}_\xi(t)$ in \eqref{deeefH}. Now, for $\xi=\frac{\epsilon}{2}$ the angle $\theta=\measuredangle(u,S^{t} u)$ is defined by

\begin{equation}\label{Angle2}
	\begin{array}{cccc}
\theta\colon& \mathbb{R}^3 & \longrightarrow & \mathbb{R} \\
& (t,\gamma,\xi) & \longmapsto & \arccos\frac{1+\gamma^2+t\xi \gamma}{\sqrt{\gamma^2+1}\sqrt{\gamma^2+1+2t\xi\gamma+(t\xi\gamma)^2}}
\end{array}
\end{equation}

\begin{figure}[h]
  \includegraphics[width=7cm,height=7cm]{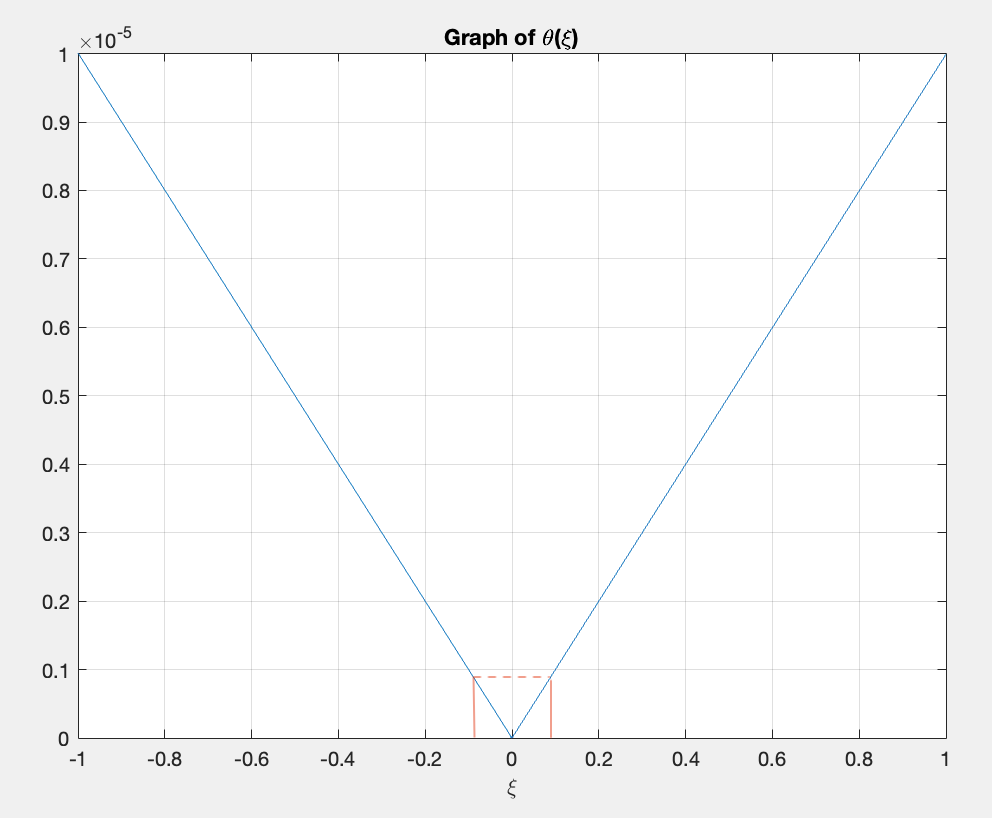}
  \caption{Graph of $\theta(\xi)$ for $\gamma=0.01$ and $\tilde\tau=0.1$. }\label{FiguraF}
\label{figure}

\end{figure}
The inequality in \eqref{second} follows directly from continuity and the map $\theta$ in \eqref{Angle2} (see Figure \ref{FiguraF}).

\end{proof}

We are now in conditions to achieve accessibility for kinetic cocycles. As highlighted at the outset of Section~\ref{esc}, our first step involves selecting a vertical cone $\mathcal{C}\gamma$. On one hand, this cone presents challenges for rotation, yet on the other hand, the original cocycle swiftly displaces vectors from the cone in time less than $\hat\tau<1/2$ (cf. Lemma~\ref{EsLem}). Once vectors exit the cone, we can engineer a suitable perturbation to redirect one direction into another (cf. Lemma~\ref{rot00} for rotation and Lemma~\ref{nota2} for selecting an appropriate time $\tilde\tau<1/2$). This manipulation relies on a map $\mathcal{A}(t)$ as defined in Definition~\ref{rot00old}, which operates effectively only outside the vertical cone $\mathcal{C}\gamma$. Overall, we are well-positioned to implement perturbations that facilitate accessibility, confined within time-1 segments of the orbit. This constitutes our next result.

	\begin{lemma}\label{rot1NEW} 
	Let $A\in L^{\infty}(M,\mathpzc{K})$ and $\epsilon>0$. There exists $\theta>0$ such that if $\w\in M$ is a $\mu$-generic point and $u,v$ belong to $\mathbb{R}^2_\w$ such that $\measuredangle(u,v)<\theta$, then there is a perturbation $B\in L^{\infty}(M,\mathpzc{K})$ of $A$ supported on $\varphi^{[0,1]}(\w)$ such that
	\begin{itemize}
	\item [(i)] $\|A(\tilde\w)-B(\tilde\w)\|<\epsilon$, for all $\tilde\w\in M$ and 
	\item [(ii)] if $U=\text{span}(u)$ and $V=\text{span}(\Phi_A(1,\w)v)$ we have $\Phi_{B}(1,\omega)U={V}$. 
	\end{itemize}
\end{lemma}

\begin{proof} 

Let be given $A\in L^{\infty}(M,\mathpzc{K})$ and $\epsilon>0$. Consider any $\hat\tau\in\left]0,\frac{1}{2}\right[$. By Lemma~\ref{EsLem} there exists $\gamma\in]0,1[$ such that given any $\w\in M$ and $v\in \mathcal{C}_\gamma\subset\mathbb{R}^2_\w$ we have $\Phi_{A}(\hat\tau,\omega)v\notin \mathcal{C}_\gamma\subset\mathbb{R}^2_{\varphi^{\hat\tau}(\w)}$.

From Lemma~\ref{nota2} there exist $\tilde\tau\in \left]0,\frac{1}{2}\right[$ and $\theta_1>0$ such that given $\w\in M$ and $u\in(\gamma,1)\cdot\mathbb{R}$ we have
\begin{equation*}\label{first2}
\max_{t\in [0,\tilde\tau]}\|H-\mathcal{A}_{\epsilon}(t)\|<\epsilon
\end{equation*}
and
\begin{equation*}\label{second2}
\measuredangle(u,R^{\tilde\tau} u)>\measuredangle(u,S^{\tilde\tau} u)=\theta_1.
\end{equation*}

Let us consider two cases:

\bigskip

\noindent\textbf{Case 1:} If $u,v\notin \mathcal{C}_\gamma$ are such that $\measuredangle(u,v)<\theta_1$, then we take $\mathcal{A}_\xi(t)$ from Definition~\ref{rot00old} and by Lemma \ref{rot00} the cocycle $B$ defined by $B=A+H$ belongs to $L^{\infty}(M,\mathpzc{K})$, where
\begin{equation}\label{Hdef900}
H(\varphi^t(\w))= \left\{
\begin{array}{ll}
      \Phi_A(t,\w)\cdot \mathcal{A}_\xi(t)\cdot [\Phi_A(t,\w)]^{-1}& \text{if}\,\,\,t\in[0,\tilde\tau] \\
      0 & \text{otherwise}
\end{array} 
\right..
\end{equation}
Moreover, for $t\in[0,\tilde\tau]$, we have $\Phi_B(t,\w)=\Phi_A(t,\w)\cdot R^t$. Therefore, we can perform the required rotation by tuning $\xi$ appropriately to get
\begin{equation}\label{RotDef}
R^{\tilde\tau} u=v
\end{equation}
From \eqref{composi} and \eqref{RotDef} we have 
\begin{equation}\label{RotDef2}
\Phi_B(\tilde\tau,\w)u=\Phi_A(\tilde\tau,\w)v.
\end{equation} 
From \eqref{Hdef900} we get 
\begin{equation}\label{RotDef3}
\Phi_B(1-\tilde\tau,\varphi^{\tilde\tau}(\w))=\Phi_A(1-\tilde\tau,\varphi^{\tilde\tau}(\w)).
\end{equation} 
Finally, \eqref{RotDef2} and \eqref{RotDef3} gives
$$\Phi_{B}(1,\omega)u=\Phi_{B}(1-\tilde\tau,\varphi^{\tilde\tau}(\w))\Phi_{B}(\tilde\tau,\omega)u=\Phi_{A}(1-\tilde\tau,\varphi^{\tilde\tau}(\w))\Phi_A(\tilde\tau,\w)v=\Phi_{A}(1,\omega)v,$$
and (ii) holds.

\bigskip
\noindent\textbf{Case 2:}

 If $u,v\in \mathcal{C}_\gamma\subset\mathbb{R}^2_\w$, then by Lemma~\ref{EsLem} given any $\w\in M$ we have $\Phi_{A}(\hat\tau,\omega)u,\Phi_{A}(\hat\tau,\omega)v\notin \mathcal{C}_\gamma\subset\mathbb{R}^2_{\varphi^{\hat\tau}(\w)}$.  By Lemma~\ref{nota1} we may choose $\theta_2>0$ such that if $\measuredangle(u,v)<\theta_2$, then
 $$\measuredangle(\Phi_{A}(\hat\tau,\omega)u,\Phi_{A}(\hat\tau,\omega)v)<\theta_1.$$ 
Notice that if there is only one vector in $\mathcal{C}_\gamma\subset\mathbb{R}^2_\w$ by Remark~\ref{rem3} there exists $\theta_0\in]0,\theta_2[$ such that if $\measuredangle(u,v)<\theta_0$, then both iterations escape from the cone.
Finally, we can reduce now to \textbf{Case 1}. The $\theta$ is defined by $\min\{\theta_0,\theta_1,\theta_2\}$.\end{proof}

\bigskip

\begin{remark}\label{consdiss} 
	Lemma~\ref{rot1NEW} also holds for the frictionless (respectively dissipative) case, because given $A\in L^{\infty}(M,\mathpzc{K}_{\,\,\Delta=0})$ (respectively $A\in L^{\infty}(M,\mathpzc{K}_{\,\,\Delta<0})$), the perturbation $B$ is obtained by considering $A+H$ where $H\in \mathbb{R}^{0\times 2}$ and $H$ is traceless, so $B \in L^{\infty}(M,\mathpzc{K}_{\,\,\Delta=0})$ (respectively $B \in L^{\infty}(M,\mathpzc{K}_{\,\,\Delta<0})$).
\end{remark}

\begin{remark}\label{re:rot}
We may extend the previous perturbations from a segment of the orbit to a flowbox $\mathcal B=\varphi^{[0,1]}(B)=\{\varphi^t(\w)\colon \w\in B,\, t\in[0,1]\}$, $B\subseteq M$, with no self-intersections, when we have attached a pair of unit vectors $u(\w),v(\w)$ for all $\w\in B$. 
\end{remark}

\bigskip

Finally, we obtain:

\begin{proposition}\label{access}
The families of kinetic cocycles in the present paper are all accessible.
\end{proposition}

\bigskip

\section{Proof of the theorems}

\subsection{Proof of the $L^\infty$ version of Theorem~\ref{tps}}

Next results (Lemmas~\ref{rot4}, \ref{rot5} and \ref{rot6}) will be stated for elements in $L^{\infty}(M,\mathpzc{K})$ but \emph{mutatis mutandis} the reader can easily see that work also for $L^{\infty}(M,\mathpzc{K}_{\,\,\Delta=0})$ and $L^{\infty}(M,\mathpzc{K}_{\,\,\Delta<0})$. Before stating the lemmas let us consider several useful sets. Let be given $A\in L^{\infty}(M,\mathpzc{K})$ over the flow $\varphi^t\colon M\rightarrow M$ and a $\varphi^t$-invariant set $\Lambda\subseteq M$.  Let $\mathcal{O}(A)=\mathcal{O}^T(A)\cup \mathcal{O}^S(A)$ be set of regular points provided by Oseledets' theorem \cite{O} where $\mathcal{O}^T(A)$ stand for the points with trivial spectrum and $\mathcal{O}^S(A)$ stand for the points with simple spectrum.  Fixing $m\in \mathbb{N}$ we say that $\Lambda$ has $m$-\emph{dominated splitting} if for $\mu$-a.e. $\w\in\Lambda$ the fiber $\mathbb{R}^2_\w$ decomposes into two $\Phi_A$-invariant one-dimensional directions $E^1_{\w}\oplus E^2_{\w}$ such that:
\begin{equation}\label{decdom}
\triangle_m(A,\w):=\frac{\left\|\Phi_A(m,\w)\Big|_{E^2_{\w}}\right\|}{\left\|\Phi_A(m,\w)\Big|_{E^1_{\w}}\right\|}\leq\frac{1}{2}.
\end{equation}

Let us fix some notation:
\begin{itemize}
\item $\mathcal{D}_m(A)$ is the set of points in $M$ whose orbit display an $m$-dominated splitting. 
\item $\Gamma_m(A)=\left\{\w\in \mathcal{O}^S(A)\colon \triangle_m(A,\w)\geq\frac{1}{2}\right\}$.
\item $\Gamma_m^\star(A)=\{\w\in\Gamma_m(A)\colon \w\text{\,is not periodic}\}$ and
\item $\Gamma_\infty(A)=\bigcap_{m\in\mathbb{N}}\Gamma_m(A)$.
\end{itemize}

The set $\Lambda$ has an $m$-\emph{uniformly hyperbolic splitting} if for $\mu$-a.e. $\w\in\Lambda$ the fiber $\mathbb{R}^2_\w$  decomposes into two $\Phi_A$-invariant one-dimensional directions $E^1_{\w}\oplus E^2{\w}$ such that:
\begin{equation}\label{hypdecdom}
\left\|\Phi_A(m,\w)^{-1}\Big|_{E^2_{\w}}\right\|\leq \frac{1}{2}\,\,\,\,\,\text{and}\,\,\,\,\,\left\|\Phi_A(m,\w)\Big|_{E^1_{\w}}\right\|\leq\frac{1}{2}.
\end{equation}

In \cite[Lemma 4.1]{BV2} was proved that $\Gamma_\infty(A)$ contains no periodic points.
Lemma~\ref{rot4} bellow is a well-know statement which in rough terms say that in the absent of a dominated splitting and along a large orbit segment we can implement rotations of a tiny angle in order to interchange two given directions. Clearly, we will apply Lemma~\ref{rot1NEW} $m$ times observing that, as the perturbations have distinct support, the concatenation of these $m$ perturbations do not add with respect to the size of the final perturbation. 

Let $\Upsilon_m(A)\subset M$ be the set of points $\w$ such that $\triangle_m(A,\w)\geq 1/2$.
Hence, $\Gamma_m(A)=\cup_{t\in \mathbb{R}}\varphi^t(\Upsilon_m(A))$.

\begin{lemma}\label{rot4} 
	Let be given $A\in L^{\infty}(M,\mathpzc{K})$ and $\epsilon>0$. There exists $m\in\mathbb{N}$ such that for $\mu$-a.e. $\w\in \Upsilon_m(A)$ there exists $B_{\w}\in L^{\infty}(M,\mathpzc{K})$ supported on the segment $\varphi^{[0,m]}(\w)$ and such that
	\begin{itemize}
	\item [(i)] $\|A-B_{\w}\|<\epsilon$ and
	\item [(ii)] $\Phi_{B_\w}(m,\w) E^2_{\w}=E^1_{\w}$.
	\end{itemize}
\end{lemma}

The second result says that interchanging directions as provided by Lemma~\ref{rot4} is effective when the task is to diminish the norm growth \emph{pointwise}.

\begin{lemma}\label{rot5} 
	Let be given $A\in L^{\infty}(M,\mathpzc{K})$, $\delta>0$ and $\epsilon>0$. If $m\in\mathbb{N}$ is sufficiently large, then there exists a measurable function $N\colon \Gamma_m^\star(A)\rightarrow\mathbb{R}$ such that for $\mu$-a.e. $\w\in\Gamma_m^\star(A)$ and every $\tau\geq N(\w)$ there exists $B_{\w}\in L^{\infty}(M,\mathpzc{K})$ supported on the segment  $\varphi^{[0,\tau]}(\w)$ and such that $\|A-B_{\w}\|<\epsilon$ and 
	\begin{equation}\label{point}
	\frac{1}{\tau}\log\|\Phi_{B_{\w}}(\tau,\w)\|\leq \delta+\frac{\lambda_1(A,\w)+\lambda_2(A,\w)}{2}.
	\end{equation}
	 Moreover, $\Phi_{B_{\w}}(t,\w)$ depends measurably on $\w$ and continuously on $t$.
\end{lemma}

\begin{proof}
For full details on the proof we refer \cite{B,BV2,Be, Be2}. The highlights are the following. We assume that $\mu(\Gamma_m^\star(A))>0$ otherwise the proof ends.  By Lemma~\ref{rot4} if $m$ is sufficiently large, then for $\w'\in\Upsilon_m(A)$ there exists a perturbation $B_{\w'}\in L^{\infty}(M,\mathpzc{K})$ supported on the segment   $\varphi^{[0,m]}(\w)$ and such that (i) $\|A-B_{\w'}\|<\epsilon$ and (ii)   $\Phi_{B_{\w'}}(m,{\w'}) E^2_{\w'}=E^1_{\w'}$. A laborious measure theoretical reasoning (see e.g. \cite[Proposition 4.2]{BV2} or \cite[Lemma 3.13]{B}) can be made to obtain a measurable function $N\colon \Gamma_m^\star(A)\rightarrow\mathbb{R}$ such that for $\mu$-a.e. $\w\in\Gamma_m^\star(A)$ and every $\tau\geq N(\w)$, there is $\ell\sim \tau/2$ such that $\w'=\varphi^{\ell}(\w)\in \Upsilon_m(A)$. The way we build the perturbation $B_\w$ is by considering:
\begin{equation}\label{mix}
\mathbb{R}^2_\w  \overset{\Phi_A(t_1,\w)}{\longrightarrow}\mathbb{R}^2_{\varphi^\ell(\w)=\w'}\overset{\Phi_{B_{\w'}}(m,\w')}{\longrightarrow}\mathbb{R}^2_{\varphi^m(\w')}\overset{\Phi_A(t_2,\w)}{\longrightarrow}\mathbb{R}^2_{\varphi^{m+t_1+t_2}(\w)},
\end{equation}
where $t_1+m+t_2=\tau$ and since for $i=1,2$, $t_i\gg m$, we have $\frac{t_i}{\tau}\sim \frac{1}{2}$. If no perturbations were introduced, we would obtain
$$\frac{1}{\tau}\log\|\Phi_{A}(\tau,\w)\|\sim \frac{\lambda_1(A,\w)+\lambda_1(A,\w)}{2}.$$ Once we perform the perturbation like in \eqref{mix} the exponential rate associated to $\lambda_1$ working on half of the way mixes with the exponential rate associated to $\lambda_2$ working on the other half. This results in the estimate in \eqref{point}.
\end{proof}

The third result, says that interchanging directions provided by Lemma~\ref{rot4} is effective when the task is to diminish the norm grow \emph{globally}. Again a laborious measure theoretical reasoning envolving Kakutani castles do the job (see e.g. \cite[Lemma 7.4]{BV2}). The extrapolation of inequality \eqref{global} from \eqref{point} is far from being an easy task. Indeed, \eqref{global} gives a global information on the Lyapunov exponents while \eqref{point} gives only pointwise information on a finite iteration.

\begin{lemma}\label{rot6} 
	Let be given $A\in L^{\infty}(M,\mathpzc{K})$, $\delta>0$ and $\epsilon>0$. There exists $m\in\mathbb{N}$ and $B\in L^{\infty}(M,\mathpzc{K})$, with  $\rho_\infty(A,B)<\epsilon$, and that equals $A$ outside the open set $\Gamma_m(A)$,  such that 
	\begin{equation}\label{global}
	\int_{\Gamma_m(A)}\lambda_1(B,\w)\,d\mu(\w)\leq \delta+\int_{\Gamma_m(A)}\frac{\lambda_1(A,\w)+\lambda_2(A,\w)}{2}\,d\mu(\w).
	\end{equation}
\end{lemma}

Since by Lemma~\ref{usc} we have that the function \eqref{uscf} is upper semicontinuous and the continuity points of upper semicontinuous functions is a residual subset we only have to check that if $A\in L^\infty(M,\mathpzc{K})$ is a continuity point of $LE(A,\Gamma_\infty(A))$ then $\lambda_1(A,\w)=\lambda_2(A,\w)$ for $\mu$-a.e. $\w\in M$. Let $\w\in M$ be an Oseledets regular point for $A$. If the two Lyapunov exponents of $A$ at $\w$ coincide, there is nothing to prove. Otherwise, $\lambda_1(A,\w)>\lambda_2(A,\w)$, $\w\notin \Gamma_\infty(A)$ and so $\w\in \mathcal{D}_m(A)$ for some $m\in \mathbb{N}$. As a conclusion we have:

\begin{theorem}\label{tps4}
		Let $\varphi^t:M\to M$ be a flow preserving the measure $\mu$. There exists a $\rho_\infty$ residual set $\mathcal{R}\subset L^\infty(M,\mathpzc{K})$, such that if $A\in \mathcal{R}$, then for $\mu$-a.e. $\omega$ either 
		\begin{itemize}
		\item there exists a single Lyapunov exponent or else
		\item the splitting along the orbit of $\w$ is dominated.
		\end{itemize}
	\end{theorem}

\subsection{Proof of Theorem~\ref{tps}}

Next result proves that a map in $L^\infty(M,\mathpzc{K})$ belongs to $C^0(M,\mathpzc{K})$ on nearly all its domain.
\begin{lemma}\label{Lusin}\textbf{Lusin-type theorem:}
Let $A\in L^\infty(M,\mathpzc{K})$ and $\epsilon>0$. There $B\in C^0(M,\mathpzc{K})$ such that $$\mu\bigl(\bigl\{\w\in M\colon A(\w)\not=B(\w)\bigr\}\bigr)<\epsilon.$$
\end{lemma}

\begin{proof}
By Lusin's theorem (see e.g \cite{Rd}) given $\alpha,\beta\in L^\infty(\mu)$ and $\epsilon>0$ there exist compact sets $M_\alpha,M_\beta\subset M$ such that $\tilde\alpha=\alpha\vert_{M_\alpha}\in C^0(M_\alpha)$, $\tilde\beta=\beta\vert_{M_\beta}\in C^0(M_\beta)$, $\mu(M_\alpha)>1-\epsilon/2$ and $\mu(M_\beta)>1-\epsilon/2$.
As the intersection of two compact subsets on a Hausdorff space is compact we define $\tilde M=M_\alpha\cap M_\beta$. By Tietze extension theorem (see e.g. \cite{Mu}) we can define $\hat\alpha\in C^0(M)$ and $\hat\beta\in C^0(M)$ as extensions, respectively, of the continuous functions $\tilde\alpha\vert_{\tilde M}$ and $\tilde\beta\vert_{\tilde M}$. Clearly, $B$ defined by 
$$B(\w)=\begin{pmatrix}0 & 1\\ -\hat\beta(\w) & -\hat\alpha(\w)\end{pmatrix}$$ belongs to $C^0(M,\mathpzc{K})$ and we have $\mu(\{\w\in M\colon A(\w)\not=B(\w)\})<\epsilon$.
\end{proof}
Lemma~\ref{Lusin} and Theorem~\ref{tps4} will imply Theorem~\ref{tps}. The way we use these results is somehow natural and next we present the highlight of the proof of a simple case (Theorem~\ref{tps2}). It is well-know that the set of continuity points of an upper semicontinuous function is a residual subset. So we will prove first that if $A\in C^0(M,\mathpzc{K})$ is a continuity point of the upper semicontinuous function $LE(\cdot, M)$ (recall Lemma~\ref{usc}), then the dichotomy in Theorem~\ref{tps2} holds. The proof is by contradiction assuming that exists a $\mu$-positive measure set $\Gamma\subset M$ without a dominated splitting and such that $LE(A,\Gamma)>0$. We then use Theorem~\ref{tps4} and construct $B_0\in L^\infty(M,\mathpzc{K})$ such that $LE(B_0,\Gamma)\sim 0$ and so $LE(B_0,\Gamma)\ll LE(A,\Gamma)$. Now, Lemma~\ref{Lusin} allows us to obtain $B\in C^0(M,\mathpzc{K})$, $\rho_\infty$-arbitrarily close to $B_0$ (thus $\rho_0$ close to $A$) and such that we still have $LE(B,\Gamma)\ll LE(A,\Gamma)$ which draw in a contradiction with the assumption the $A$ was a continuity point of $LE(\cdot, M)$. Here, $\ll$ means that we can cause a discontinuity on the map $LE$. 

\subsection{Proof of Theorem~\ref{tps2} and Theorem~\ref{tps3}}

To prove Theorem~\ref{tps2} we notice that our perturbations keep us inside $L^{\infty}(M,\mathpzc{K}_{\,\,\Delta=0})$. As $\lambda_1(A,\w)=-\lambda_2(A,\w)$ for $\mu$-a.e. $\w\in M$ the version of Lemma~\ref{rot6} will be as follows: Let be given $A\in L^{\infty}(M,\mathpzc{K}_{\,\,\Delta=0})$, $\delta>0$ and $\epsilon>0$. There exists $m\in\mathbb{N}$ and $B\in L^{\infty}(M,\mathpzc{K}_{\,\,\Delta=0})$ with $\|A-B_{\w}\|<\epsilon$ that equals $A$ outside the open set $\Gamma_m(A)$ and such that $\int_{\Gamma_m(A)}\lambda_1(B,\w)\,d\mu(\w)\leq \delta$. To get the $C^0$ version of our result we use a similar reasoning like in Lemma~\ref{Lusin}. 
	
To prove Theorem~\ref{tps3} we use again the perturbations developed in \S\ref{toolbox} to rotate solutions and remain inside $L^{\infty}(M,\mathpzc{K}_{\,\,\Delta<0})$. Now, the version of Lemma~\ref{rot6} will be as follows: Let be given $A\in L^{\infty}(M,\mathpzc{K}_{\,\,\Delta<0})$, $\delta>0$ and $\epsilon>0$. There exists $m\in\mathbb{N}$ and $B\in L^{\infty}(M,\mathpzc{K}_{\,\,\Delta<0})$ that equals $A$ outside the open set $\Gamma_m(A)$, with $\rho_\infty(A,B)<\epsilon$, and such that \eqref{global} holds. Again on $\Gamma_\infty(A)$ we get the first item of Theorem~\ref{tps3} i.e. the solution of almost every equation \eqref{SOLH} is stable or equivalently, for $\mu$-a.e. the exists a single negative Lyapunov exponent. Otherwise we have $\mu$-a.e. $\w\in \mathcal{D}_m(A)$ for some $m\in \mathbb{N}$. From Oseledets' theorem and Liouville's formula we get:
\begin{equation*}
\lambda_1(A,\w)+\lambda_2(A,\w)=\lim\frac{1}{t}\log\left|\det \Phi_A(t,\w)\right|=\lim\frac{1}{t}\int_0^t \text{Trace}(A(\varphi^s(\w)))\,ds<0.
\end{equation*}
Now, we have two possible situations: the solution of every associated equation \eqref{SOLH} is stable (dominated splitting with two negative Lyapunov exponents) or else the solution of every associated equation \eqref{SOLH} displays a dominated splitting with two Lyapunov exponents of different signs. Hence, for any initial condition in the complement of the Oseledets direction associated to the negative Lyapunov exponent the solution of every associated equation \eqref{SOLH} is uniformly unstable which is precisely the last item of Theorem~\ref{tps3}.

\subsection{Theorem~\ref{tps2} applied to the 1-d Schr\"odinger case}\label{Sc}

Finally, we present traceless continuous-time linear cocycles with a somewhat different look, namely by considering the one-di\-men\-sional Schr\"odin\-ger operator on $L^2(\mathbb{R})$ and with an $C^0$ potential $Q\colon M\rightarrow \mathbb{R}$ given by:
	\begin{equation}\label{Sch}
	\begin{array}{crcl}
H_\w\colon &L^2(\mathbb{R}) & \longrightarrow & L^2(\mathbb{R}) \\& \psi & \longmapsto &  \left[-\frac{d^2}{dt^2}+Q(\varphi^t(\w))\right]\psi
\end{array}
	\end{equation}
	In particular we like to describe the Lyapunov spectrum of the time-independent Schr\"odinger equation
	\begin{equation}\label{Sch2}
	H_\w\psi=E\psi,
	\end{equation}
	where $E\in \mathbb{R}$ is a given energy. Putting together \eqref{Sch} and \eqref{Sch2} we deduce a kinetic cocycle as in \eqref{damp2} but with $\alpha(\w)=0$ and $\beta(\w)=E-Q(\w)$ for all $\w\in M$. We fix the energy $E$ and focus on the continuous-time linear cocycle
	\begin{equation*}\label{damp22}
\begin{array}{cccc}
A_E\colon& M & \longrightarrow &  \mathbb{R}^{2\times2} \\
& \w & \longmapsto & \begin{pmatrix}0 & 1\\ -E+Q(\omega) & 0\end{pmatrix}
\end{array}
\end{equation*}
called 1-d Schr\"odinger LDS with potential $Q$.

Let $C^0(M,\mathbb{R})$ stand for the set of continuous potentials. As a direct consequence of Theorem~\ref{tps2} we have:
	
		\begin{maintheorem}\label{tps5}
		Let $\varphi^t:M\to M$ be a flow preserving the measure $\mu$ and fix the energy $E$. There exists a $\rho_0$ residual set $\mathcal{R}\subset C^0(M,\mathbb{R})$, such that if $Q\in \mathcal{R}$, then for the 1-d Schr\"odinger continuous-time linear cocycle with energy $E$ and potential $Q$ either 
		\begin{itemize}
		\item  all Lyapunov exponents vanishes or else
		\item the splitting along the orbit of $\w$ is uniformly hyperbolic.
		\end{itemize}
	\end{maintheorem}

\vspace{1cm}
	
\textbf{Acknowledgements:} \small {MB was partially supported by CMUP, which is financed by national funds through FCT-Funda\c c\~ao para a Ci\^encia e a Tecnologia, I.P., under the project with reference UIDB/00144/2020. MB was also partially supported by
the project \emph{Means and Extremes in Dynamical Systems} PTDC/MAT-PUR/4048/2021. HV was partially supported by FCT - `Funda\c{c}\~ao para a Ci\^encia e a Tecnologia', through Centro de Matem\'atica e Aplica\c{c}\~oes (CMA-UBI), Universidade da Beira Interior, project UIDB/MAT/00212/2020. }

\end{document}